\documentclass[11pt,leqno]{article}

\usepackage{amsmath,amsfonts,amscd,amssymb,theorem}

\long\def\comment#1\endcomment{}

\comment
\endcomment

\makeatletter
\begingroup
\gdef\th@dotted{\normalfont\itshape
  \def\@begintheorem##1##2{%
        \item[\hskip\labelsep \theorem@headerfont ##1\ ##2.]}%
\def\@opargbegintheorem##1##2##3{%
   \item[\hskip\labelsep \theorem@headerfont ##1\ ##2\ (##3).]}}
\endgroup
\makeatother

\theoremstyle{dotted}

\newtheorem{theorem}{Theorem}[section]
\newtheorem{lemma}[theorem]{Lemma}

\makeatletter
\begingroup
\gdef\th@upshape{\normalfont
  \def\@begintheorem##1##2{%
        \item[\hskip\labelsep \theorem@headerfont ##1\ ##2.]}%
\def\@opargbegintheorem##1##2##3{%
   \item[\hskip\labelsep \theorem@headerfont ##1\ ##2\ (##3).]}}
\endgroup
\makeatother

\theoremstyle{upshape}

\newtheorem{defn}[theorem]{Definition}
\newtheorem{remark}[theorem]{Remark}
\newtheorem{exa}[theorem]{Example}

\makeatletter
\renewcommand{\subsection}{\@startsection{subsection}{2}{0pt}{-3ex
plus -1ex minus -0.2ex}{-2mm plus -0pt minus
-2pt}{\normalfont\bfseries}} 
\renewcommand{\subsubsection}{\@startsection{subsubsection}{3}{0pt}{-3ex
plus -1ex minus -0.2ex}{-2mm plus -0pt minus
-2pt}{\normalfont\bfseries}} 
\makeatother

\makeatletter
\@addtoreset{equation}{section}
\makeatother

\newcommand{\cntrct}                
{\hspace{2pt}\raisebox{1pt}{\text{$\lrcorner$}}\hspace{2pt}}

\newcommand{\proof}[1][Proof.]{\smallskip\noindent{\em #1}}
\def\endproof{\hfill\ensuremath{\square}\par\medskip}

\def\eqref#1{\thetag{\ref{#1}}}

\let\latexref=\ref
\def\ref#1{{\normalfont{\latexref{#1}}}}

\newcommand{\wt}{\widetilde}

\newcommand{\idot}{{\:\raisebox{1pt}{\text{\circle*{1.5}}}}}
\newcommand{\hdot}{{\:\raisebox{3pt}{\text{\circle*{1.5}}}}}
\renewcommand{\phi}{\varphi}

\def\dlim_#1{{\displaystyle\lim_{#1}}^\hdot}

\newcommand{\Hom}{\operatorname{Hom}}

\newcommand{\RHom}{\operatorname{RHom}}
\newcommand{\Tor}{\operatorname{Tor}}

\newcommand{\id}{\operatorname{\sf id}}

\newcommand{\rk}{\operatorname{\sf rk}}
\newcommand{\Rk}{\operatorname{\it rk}}

\newcommand{\Tot}{\operatorname{Tot}}

\newcommand{\D}{{\cal D}}
\newcommand{\C}{{\cal C}}

\newcommand{\Q}{{\cal Q}}

\newcommand{\K}{\mathcal{K}}

\newcommand{\Sets}{\operatorname{Sets}}
\newcommand{\Cat}{\operatorname{Cat}}
\newcommand{\Top}{\operatorname{Top}}
\newcommand{\Aut}{{\operatorname{Aut}}}

\newcommand{\hhom}{\operatorname{\mathcal{H}{\it om}}}

\newcommand{\Add}{\operatorname{\sf Add}}

\newcommand{\amod}{{\text{\rm -mod}}}
\newcommand{\ppt}{{\sf pt}}

\newcommand{\lotimes}{\overset{\sf\scriptscriptstyle L}{\otimes}}

\newcommand{\cchar}{\operatorname{\sf char}}

\newcommand{\Mat}{{\sf Mat}}
\newcommand{\mat}{\mathcal{M}{\it at}}
\newcommand{\Vect}{{\sf Vect}}
\newcommand{\vect}{\mathcal{V}{\it ect}}

\newcommand{\Dd}{\operatorname{\sf D}}

\newcommand{\Z}{{\mathbb Z}}

\newcommand{\QQ}{{\mathbb Q}}

\newcommand{\Fun}{\operatorname{Fun}}
\newcommand{\G}{\operatorname{Gr}}

\newcommand{\E}{\mathcal{E}}

\newcommand{\N}{\mathcal{N}}
\newcommand{\Nn}{\operatorname{\sf N}}

\newcommand{\kk}{\mathbb{K}}
\newcommand{\ii}{\overline{i}}

\newcommand{\Ss}{\mathcal{S}}
\newcommand{\X}{\mathcal{X}}

\title{$K$-theory as an Eilenberg-Maclane spectrum}

\author{D. Kaledin\thanks{Partially supported by RScF, grant number
    14-21-00053, by AG Laboratory SU-HSE, RF government grant,
    ag.11.G34.31.0023, and by the Dynasty Foundation award.}}

\date{To Sasha Merkurjev, on the occasion of his 60th anniversary}

\begin{document}

\maketitle

\tableofcontents

\section*{Introduction.}

Various homology and cohomology theories in algebra or algebraic
geometry usually take as input a ring $A$ or an algebraic variety
$X$, and produce as output a certain chain complex; the homology
groups of this chain complex are by definition the homology or
cohomology groups of $A$ or $X$. Higher algebraic $K$-groups are
very different in this respect -- by definition, the groups
$K_\idot(A)$ are homotopy groups of a certain spectum $\K(A)$. Were
it possible to represent $K_\idot(A)$ as homology groups of a chain
complex, one would be able to study it by means of the
well-developed and powerful machinery of homological
algebra. However, this is not possible: the spectrum $\K(A)$ is
almost never a spectrum of the Eilenberg-Maclane type.

If one wishes to turn $\K(A)$ into an Eilenberg-Maclane spectrum,
one needs to complete it or to localize it in a certain set of
primes. The cheapest way to do it is of course to localize in {\em
  all} primes -- rationally, every spectrum is an Eilenberg-Maclane
spectrum, and the difference between spectra and complexes
disappears. The groups $K_\idot(A) \otimes \QQ$ are then the
primitive elements in the homology groups
$H_\idot(BGL_{\infty}(A),\QQ)$, and this allows for some
computations using homological methods. In particular, $K_\idot(A)
\otimes \mathbb{Q}$ has been computed by Borel when $A$ is a number
field, and the relative $K$-groups $K_\idot(A,I) \otimes \QQ$
of a $\QQ$-algebra $A$ with respect to a nilpotent ideal $I
\subset A$ have been computed in full generality by Goodwillie
\cite{good}.

However, there is at least one other situation when $\K(A)$ becomes
an Eilenberg-Maclane spectrum after localization. Namely, if $A$ is
a finite field $k$ of characteristic $p$, then by a famous result of
Quillen \cite{qui}, the localization $\K^{(p)}(A)$ of the spectrum
$\K(A)$ at $p$ is the Eilenberg-Maclane spectrum $H(\Z_{(p)})$
corresponding to the ring $\Z_{(p)}$. Moreover, if $A$ is an algebra
over $k$, then $\K(A)$ is a module spectum over $\K(k)$ by a result
of Gillet \cite{gil}. Then $\K^{(p)}(A)$ is a module spectrum over
$H(\Z_{(p)})$, thus an Eilenberg-Maclane spectrum corresponding to a
chain complex $K_\idot^{(p)}(A)$ of $\Z_{(p)}$-modules. More
generally, if we have a $k$-linear exact or Waldhausen category
$\C$, the $p$-localization $\K^{(p)}(\C)$ of the $K$-theory spectrum
$\K(\C)$ is also an Eilenberg-Maclane spectrum corresponding to a
chain complex $K^{(p)}_\idot(\C)$.

Moreover, if we have a nilpotent ideal $I \subset A$ in a
$k$-algebra $A$, then the relative $K$-theory spectrum $\K(A,I)$ is
automatically $p$-local. Thus $\K(A,I) \cong \K^{(p)}(A,I)$ is an
Eilenberg-Maclane spectrum ``as is'', without further modifications.

Unfortunately, unlike in the rational case, the construction of the
chain complex $K_\idot^{(p)}(\C)$ is very indirect and uncanonical,
so it does not help much in practical computations. One clear
deficiency is insufficient functoriality of the construction that
makes it difficult to study its behaviour in families. Namely, a
convenient axiomatization of the notion of a family of categories
indexed by small category $\C$ is the notion of a cofibered category
$\C'/\C$ introduced in \cite{SGA}. This is basicaly a functor
$\pi:\C' \to \C$ satisfying some conditions; the conditions insure
that for every morphism $f:c \to c'$ in $\C'$, one has a natural
transition functor $f_!:\pi^{-1}(c) \to \pi^{-1}(c')$ between fibers
of the cofibration $\pi$. Cofibration also behave nicely with
respect to pullbacks -- for any cofibered category $\C'/\C$ and any
functor $\gamma:\C_1 \to \C$, we have the induced cofibration
$\gamma^*\C' \to \C_1$. Within the context of algebraic $K$-theory,
one would like to start with a cofibration $\pi:\C' \to \C$ whose
fibers $\pi^{-1}(c)$, $c \in \C$, are $k$-linear additive
categories, or maybe $k$-linear exact or Waldhausen categories, and
one would like to pack the individual complexes
$K^{(p)}_\idot(\pi^{-1}(c))$ into a single object
$$
K^{(p)}(\C'/\C) \in \D(\C,\Z_{(p)})
$$
in the derived category $\D(\C,\Z_{(p)})$ of the category of
functors from $\C$ to $\Z_{(p)}$-modules. One would also like this
construction to be functorial with respect to pullbacks, so that for
any functor $\gamma:\C_1 \to \C$, one has a natural isomorphism
$$
\gamma^*K^{(p)}(\C'/\C) \cong K^{(p)}(\gamma^*\C'/\C_1).
$$
In order to achieve this by the usual methods, one has to construct
the chain complex $K^{(p)}_\idot(\C)$ in such a way that it is
exactly functorial in $\C$. This is probably possible but extremely
painful.

\medskip

The goal of this paper, then, is to present an alternative very
simple construction of the objects $K^{(p)}(\C'/\C) \in
\D(\C,\Z_{(p)})$ that only uses direct homological methods, without
any need to even introduce the notion of a ring spectrum.  The only
thing we need to set up the construction is a commutative ring $k$
and a localization $R$ of the ring $\Z$ in a set of primes $S$ such
that $K_i(k) \otimes R = 0$ for $i \geq 1$, and $K_0(k) \otimes R
\cong R$. Starting from these data, we produce a family of objects
$K^R(\C'/\C) \in \D(\C,R)$ with the properties listed above, and
such that if $\C$ is the point category $\ppt$, then $K^R(\C'/\ppt)$
is naturally identified with the $K$-theory spectrum $\K(\C')$
localized in $S$.

\medskip

Although the only example we have in mind is $k = \mathbb{F}_q$, $R
= \Z_{(p)}$, we formulate things in bigger generality to emphasize
the essential ingredients of the construction. We do not need any
information on how the isomorphism $K_0(k) \otimes R \cong R$ comes
about, nor on why the higher $K$-groups vanish. As our entry point
to algebraic $K$-theory, we use the formalism of Waldhausen
categories, since it is the most general one available. However,
were one to wish to use, for example, Quillen's $Q$-construction,
everything would work with minimal modifications.

Essentially, our approach is modeled on the approach to Topological
Hochschild Homology pioneered by M. Jibladze and T. Pirashvili
\cite{JP}. The construction itself is quite elementary. The
underlying idea is also rather transparent and would work in much
larger generality, but at the cost of much more technology to make
things precise. Thus we have decided to present both the idea and
its implementation but to keep them separate. In
Section~\ref{heu.sec}, we present the general idea of the
construction, without making any mathematical statements precise
enough to be proved. The rest of the paper is completely indepedent
of Section~\ref{heu.sec}. A rather long Section~\ref{pre.sec}
contains the list of preliminaries; everything is elementary and
well-known, but we need to recall these things to set up the
notation and make the paper self-contained. A short explanation of
what is needed and why is contained in the end of
Section~\ref{heu.sec}. Then Section~\ref{sta.sec} gives the exact
statement of our main result, Theorem~\ref{main}, and
Section~\ref{proof.sec} contains its proof.

\section{Heuristics.}\label{heu.sec}

Assume given a commutative ring $R$, and let $M(R)$ be the category
of finitely generated free $R$-modules. It will be useful to
interpret $M(R)$ as the category of matrices: objects are finite
sets $S$, morphisms from $S$ to $S'$ are $R$-valued matrices of size
$S \times S'$.

Every $R$-module $M$ defines a $R$-linear additive functor $\wt{M}$
from $M(R)$ to the category of $R$-modules by setting
\begin{equation}\label{wt.m}
\wt{M}(M_1) = \Hom_R(M_1^*,M)
\end{equation}
for any $M_1 \in M(R)$, where we denote by $M_1^* = \Hom_R(M_1,R)$
the dual $R$-modules. This gives an equivalence of categories
between the category $R\amod$ of $R$-modules, and the category of
$R$-linear additive functors from $M(R)$ to $R\amod$.

Let us now make the following observation. If we forget the
$R$-module structure on $M$ and treat it as a set, we of course lose
information. However, if we do it pointwise with the functor
$\wt{M}$, we can still recover the original $R$-module $M$. Namely,
denote by $\Fun(M(R),R)$ the category of all functors from $M(R)$ to
$R\amod$, without any additivity or linearity conditions, and
consider the functor $R\amod \to \Fun(M(R),R)$ that sends $M$ to
$\wt{M}$. Then it has a left-adjoint functor
$$
\Add_R:\Fun(M(R),R) \to R\amod,
$$
and for any $M \in R\amod$, we have
\begin{equation}\label{add.m}
M \cong \Add_R(R[\wt{M}]),
\end{equation}
where $R[\wt{M}] \in \Fun(M(R),R)$ sends $M_1 \in M(R)$ to the free
$R$-module $R[\wt{M}(M_1)]$ generated by $\wt{M}(M_1)$. Indeed, by
adjunction, $\Add_R$ commutes with colimits, so it suffices to check
\eqref{add.m} for a finitely generated free $R$-module $M$; but then
$R[\wt{M}]$ is a representable functor, and \eqref{add.m} follows from
the Yoneda Lemma.

The functor $\Add_R$ also has a version with coefficients. If we
have an $R$-algebra $R'$, then for any $R'$-module $M$, the functor
$\wt{M}$ defined by \eqref{wt.m} is naturally a functor from $M(R)$
to $R'\amod$. Then by adjunction, we can define the functor
$$
\Add_{R,R'}:\Fun(M(R),R') \to R'\amod,
$$
and we have an isomorphism
\begin{equation}\label{add.r}
\Add_{R,R'}(R'[\wt{M}]) \cong M \otimes_R R'
\end{equation}
for any flat $R$-module $M$.

\medskip

What we want to do now is to obtain a homotopical version of the
construction above. We thus replace sets with topological spaces. An
abelian group structure on a set becomes an infinite loop space
structure on a topological space; this is conveniently encoded by a
special $\Gamma$-space of G. Segal \cite{segal}. Abelian groups
become connective spectra. Rings should become ring spectra. As far
as I know, Segal machine does not extend directly to ring spectra --
to describe ring spectra, one has to use more complicated machinery
such as ``functors with smash products'', or an elaboration on them,
ring objects in the category of symmetric spectra of
\cite{sym}. However, in practice, if we are given a connective
spectrum $\X$ represented by an infinite loop space $X$, then a ring
spectrum structure on $\X$ gives rise to a multiplication map $\mu:X
\times X \to X$, and in ideal situation, this is sufficiently
associative and distributive to define a matrix category $\Mat(X)$
analogous to $M(R)$. This should be a small category enriched over
topological spaces. Its objects are finite sets $S$, and the space
of morphisms from $S$ to $S'$ is the space $X^{S \times S'}$ of
$X$-valued matrices of size $S \times S'$, with compositions induced
by the multiplication map $\mu:X \times X \to X$.

Ideal situations seem to be rare (the only example that comes to
mind readily is a simplicial ring treated as an Eilenberg-Maclane
ring spectrum). However, one might relax the conditions
slightly. Namely, in practice, infinite loop spaces and special
$\Gamma$-spaces often appear as geometric realizations of monoidal
categories. The simplest example of this is the sphere spectrum
$\Ss$. One start with the groupoid $\overline{\Gamma}$ of finite
sets and isomorphisms between them, one treats it as a monoidal
category with respect to the disjoint union operation, and one
produces a special $\Gamma$-space with underlying topological space
$|\overline{\Gamma}|$, the geometric realization of the nerve of the
category $\overline{\Gamma}$. Then by Barratt-Quillen Theorem, the
corresponding spectrum is exactly $\Ss$.

The sphere spectrum is of course a ring spectrum, and the
multiplication operation $\mu$ also has a categorical origin: it is
induced by the cartesian product functor $\overline{\Gamma} \times
\overline{\Gamma} \to \overline{\Gamma}$. This functor is not
associative or commutative on the nose, but it is associative and
commutative up to canonical isomorphisms. The hypothetical matrix
category $\Mat(|\overline{\Gamma}|)$ is then easily constructed as
the geometric realization $|\Q\Gamma|$ of a strictification of a
$2$-category $\Q\Gamma$ whose objects are finite sets $S$, and whose
category $\Q\Gamma(S,S')$ of morphisms from $S$ to $S'$ is the
groupoid $\overline{\Gamma}^{S \times S'}$. Equivalently,
$\Q\Gamma(S,S')$ is the category of diagrams
\begin{equation}\label{dd}
\begin{CD}
S @<{l}<< \wt{S} @>{r}>> S'
\end{CD}
\end{equation}
of finite sets, and isomorphisms between these
diagrams. Compositions are obtained by taking pullbacks.

Any spectrum is canonically a module spectrum over $\Ss$. So, in
line with the additivization yoga described above, we expect to be
able to start with a connective spectrum $\X$ corresponding to an
infinite loop space $X$, produce a functor $X_\idot$ from $|\Q\Gamma|$
to topological spaces sending a finite set $S$ to $X^S$, and then
recover the infinite loop space structure on $X$ from the functor
$X_\idot$.

This is exactly what happens -- and in fact, we do not need the
whole $2$-category $\Q\Gamma$, it suffices to restrict our attention
to the subcategory in $\Q\Gamma$ spanned by diagrams \eqref{dd} with
injective map $l$. Since such diagrams have no non-trivial
automorphisms, this subcategory is actually a $1$-category. It is
equivalent to the category $\Gamma_+$ of pointed finite sets. Then
restricting $X_\idot$ to $\Gamma_+$ produced a functor from
$\Gamma_+$ to topological spaces, that is, precisely a
$\Gamma$-space in the sense of Segal. This $\Gamma$-space is
automatically special, and one recovers the infinite loop space
structure on $X$ by applying the Segal machine.

It is also instructive to do the versions with coefficients, with
$R$ being the sphere spectrum, and $R'$ being the Eilenberg-Maclane
ring spectrum $H(A)$ corresponding to a ring $A$. Then module
spectra over $H(A)$ are just complexes of $A$-modules, forming the
derived category $\D(A)$ of the category $A\amod$, and functors from
$\Gamma_+$ to $H(A)$-module spectra are complexes in the category
$\Fun(\Gamma_+,A)$ of functors from $\Gamma_+$ to $A\amod$, forming
the derived category $\D(\Gamma_+,A)$ of the abelian category
$\Fun(\Gamma_+,A)$. One has a tautological functor from $A\amod$ to
$\Fun(\Gamma_+,A)$ sending an $A$-module $M$ to $\wt{M} \in
\Fun(\Gamma_+,A)$ given by $\wt{M}(S) = M[\overline{S}]$, where
$\overline{S} \subset S$ is the complement to the distinguished
element $o \in S$. This has a
left-adjoint functor
$$
\Add:\Fun(\Gamma_+,A) \to A\amod,
$$
with its derived functor $L^\hdot\Add:\D(\Gamma_+,A) \to \D(A)$.
The role of the free $A$-module $A[S]$ generated by a set $S$ is
played by the singular chain complex $C_\idot(X,A)$ of a topological
space, and we expect to start with a special $\Gamma$-space
$X_+:\Gamma_+ \to \Top$, and obtain an analog of \eqref{add.r},
namely, an isomorphism
$$
L^\hdot\Add(C_\idot(X_+,A)) \cong H_\idot(\X,A),
$$
where $H_\idot(\X,A)$ are the homology groups of the spectum $\X$
corresponding to $X_+$ with coefficients in $A$ (that is, homotopy
groups of the product $\X \wedge A$).

Such an isomorphism indeed exists; we recall a precise statement
below in Lemma~\ref{loop}.

Moreover, we can be more faithful to the original construction and
avoid restricting to $\Gamma_+ \subset \Q\Gamma$. This entails a
technical difficulty, since one has to explain what is a functor
from the $2$-category $\Q\Gamma$ to complexes of $A$-modules, and
define the corresponding derived category $\D(\Q\Gamma,A)$. It can
be done in several equivalent ways, see e.g. \cite[Section
  3.1]{mackey}, and by \cite[Lemma 3.4(i)]{mackey}, the answer
remains the same -- we still recover the homology groups
$H_\idot(\X,A)$.

\medskip

Now, the point of the present paper is the following. The $K$-theory
spectrum $\K(k)$ of a commutative ring $k$ also comes from a
monoidal category, namely, the groupoid $BGL(k)$ of finitely
generated projective $k$-modules and isomorphisms between
them. Moreover, the ring structure on $\K(k)$ also has categorical
origin -- it comes from the tensor product functor $BGL(k) \times
BGL(k) \to BGL(k)$. And if we have some $k$-linear Waldhausen
category $\C$, then the infinite loop space corresponding to the
$K$-theory spectrum $\K(\C)$ is the realization of the nerve of a
category $S\C$ on which $BGL(k)$ acts. Therefore one can construct a
$2$-category $\Mat(k)$ of matrices over $BGL(k)$, and $\C$ defines a
$2$-functor $\Vect(S\C):\Mat(k) \to \Cat$ to the $2$-category $\Cat$
of small categories. At this point, we can forget all about ring
spectra and module spectra, define an additivization functor
$$
\Add:\D(\Mat(k),R) \to \D(R),
$$
and use an analog of \eqref{add.r} to recover if not $\K(\C)$ then
at least $\K(\C) \wedge_{\K(k)} H(R)$, where $H(R)$ is the
Eilenberg-Maclane spectrum corresponding to $R$. This is good
enough: if $R$ is the localization of $\Z$ in a set of primes $S$
such that $\K(k)$ localized in $S$ is $H(R)$, then $\K(\C)
\wedge_{\K(k)} H(R)$ is the localization of $\K(\C)$ in $S$.

The implementation of the idea sketched above requires some
preliminaries. Here is a list. Subsection~\ref{ho.subs} discusses
functor categories, their derived categories and the like; it is
there mostly to fix notation. Subsection~\ref{gro.subs} recalls the
basics of the Grothendieck construction of
\cite{SGA}. Subsection~\ref{bc.subs} contains some related
homological facts. Subsection~\ref{del.subs} recalls some standard
facts about simplicial sets and nerves of
$2$-categories. Subsection~\ref{2.subs} discusses $2$-categories and
their nerves. Subsection~\ref{ho.2.subs} constructs the derived
category $\D(\C,R)$ of functors from a small $2$-category $\C$ to
the category of modules over a ring $R$; this material is slightly
non-standard, and we have even included one statement with a
proof. We use an approach based on nerves, since it is cleaner and
does not require any strictification of $2$-categories. Then we
introduce the $2$-categories we will need: Subsection~\ref{fin.subs}
is concerned with the $2$-category $\Q\Gamma$ and its subcategory
$\Gamma_+ \subset \Q\Gamma$, while Subsection~\ref{mat.subs}
explains the matrix $2$-categories $\Mat(k)$ and the $2$-functors
$\Vect(\C)$. Finally, Subsection~\ref{rel.subs} explains how the
matrix and vector categories are defined in families (that is, in
the relative setting, with respect to a cofibration in the sense of
\cite{SGA}).

Having finished with the preliminaries, we turn to our
results. Section~\ref{sta.sec} contains a brief recollection on
$K$-theory, and then the statement of the main result,
Theorem~\ref{main}. Since we do not introduce ring spectra, we
cannot really state that we prove a spectral analog of
\eqref{add.r}. Instead, we construct directly a map $\K(\C) \to \K$
to a certain Eilenberg-Maclane spectrum $\K$, and we prove that the
map becomes an isomorphism after the appropriate localization. The
actual proof is contained in Section~\ref{proof.sec}.

\section{Preliminaries.}\label{pre.sec}

\subsection{Homology of small categories.}\label{ho.subs}

For any two objects $c,c' \in \C$ in a category $\C$, we will denote
by $\C(c,c')$ the set of maps from $c$ to $c'$. For any category
$\C$, we will denote by $\C^o$ the opposite category, so that
$\C(c,c') = \C^o(c',c)$, $c,c' \in \C$. For any functor $\pi:\C_1
\to \C_2$, we denote by $\pi^o:\C_1^o \to \C_2^o$ the induced
functor between the opposite categories.

For any small category $\C$ and ring $R$, we will denote by
$\Fun(\C,R)$ the abelian category of functors from $\C$ to the
category $R\amod$ of left $R$-modules, and we will denote by
$\D(\C,R)$ its derived category. The triangulated category
$\D(\C,R)$ has a standard $t$-structure in the sense of \cite{BBD}
whose heart is $\Fun(\C,R)$. For any object $c \in \C$, we will
denote by $R_c \in \Fun(\C,R)$ the representable functor given by
\begin{equation}\label{r.c}
R_c(c')=R[\C(c,c')],
\end{equation}
where for any set $S$, we denote by $R[S]$ the free $R$-module
spanned by $S$. Every object $E \in \D(\C,R)$
defines a functor $\D(E):\C \to \D(R)$ from $\C$ to the derived
category $\D(R)$ of the category $R\amod$, and by adjunction, we
have a quasiisomorphism
\begin{equation}\label{r.c.re}
\D(E)(c) \cong \RHom^\hdot(R_c,E)
\end{equation}
for any object $c \in \C$ (we will abuse notation by writing $E(c)$
instead of $\D(E)(c)$). Any functor $\gamma:\C \to \C'$ between
small categories induces an exact pullback functor
$\gamma^*:\Fun(\C',R) \to \Fun(\C,R)$ and its adjoints, the left and
right Kan extension functors $\gamma_!,\gamma_*:\Fun(\C,R) \to
\Fun(\C',R)$. The derived functors
$L^\hdot\gamma_!,R^\hdot\gamma_*:\D(\C,R) \to \D(C',R)$ are left
resp. right-adjoint to the pullback functor $\gamma^*:\D(C',R) \to
\D(\C,R)$. The homology resp. cohomology of a small category $\C$
with coefficients in a functor $E \in \Fun(\C,R)$ are given by
$$
H_i(\C,E) = L^i\tau_! E, \qquad H^i(\C,E) =
R^i\tau_*E, \qquad i \geq 0,
$$
where $\tau:\C \to \ppt$ is the tautological projection to the point
category $\ppt$.

Assume that the ring $R$ is commutative. Then for any $E \in
\Fun(\C,R)$, $T \in \Fun(\C^o,R)$, the {\em tensor product} $E
\otimes_\C T$ is the cokernel of the natural map
$$
\begin{CD}
\displaystyle\bigoplus_{f:c \to c'}E(c) \otimes_R T(c') @>{E(f)
  \otimes \id - \id \otimes T(f)}>> \displaystyle\bigoplus_{c \in
  \C}E(c) \otimes_R T(c).
\end{CD}
$$
Sending $E$ to $E \otimes_\C T$ gives a right-exact functor from
$\Fun(\C,R)$ to $R\amod$; we denote its derived functors by
$\Tor_i^\C(E,T)$, $i \geq 1$, and we denote by $E \lotimes T$ the
derived tensor product. If $T(c)$ is a free $R$-module for any $c
\in \C$, then $- \otimes_\C T$ is left-adjoint to an exact functor
$\hhom(T,-):R\amod \to \Fun(\C,R)$ given by
\begin{equation}\label{hhom}
\hhom(T,E)(c) = \Hom(T(c),E), \qquad c \in \C,E \in R\amod.
\end{equation}
Being exact, $\hhom(T,-)$ induces a functor from $\D(R)$ to
$\D(\C,R)$; this functor is right-adjoint to the derived tensor
product functor $- \lotimes_\C T$. For example, if $T = R$ is the
constant functor with value $R$, then we have
$$
H_\idot(\C,E) \cong \Tor_\idot^\C(E,\Z)
$$
for any $E \in \Fun(\C,R)$.

\subsection{Grothendieck construction.}\label{gro.subs}

A morphism $f:c \to c'$ in a category $\C'$ is called {\em
  cartesian} with respect to a functor $\pi:\C' \to \C$ if any
morphism $f_1:c_1 \to c'$ in $\C'$ such that $\pi(f) = \pi(f_1)$
factors uniquely as $f_1 = f \circ g$ for some $g:c_1 \to c$. A
functor $\pi:\C' \to \C$ is a {\em prefibration} if for any morphism
$f:c \to c'$ in $\C$ and object $c_1' \in \C'$ with $\pi(c_1') =
c'$, there exists a cartesian map $f_1:c_1 \to c_1'$ in $\C'$ with
$\pi(f_1)=f$. A prefibration is a {\em fibration} if the composition
of two cartesian maps is cartesian. A functor $F:\C' \to \C''$
between two fibrations $\C',\C''/\C$ is {\em cartesian} if it
commutes with projections to $\C$ and sends cartesian maps to
cartesian maps. For any fibration $\C' \to \C$, a subcategory $\C'_0
\subset \C'$ is a {\em subfibration} if the induced functor $\C'_0
\to \C$ is a fibration, and the embedding functor $\C'_0 \to \C'$ is
cartesian over $\C$.

A fibration $\pi:\C' \to \C$ is called {\em discrete} if its fibers
$\pi_c = \pi^{-1}(c)$, $c \in \C$ are discrete categories. For
example, for any $c \in \C$, let $\C/c$ be the category of objects
$c' \in \C$ equipped with a map $c' \to c$. Then the forgetful
functor $\phi:\C/c \to \C$ sending $c' \to c$ to $c'$ is a discrete
fibration, with fibers $\phi_{c'} = \C(c',c)$, $c' \in \C$.

For any functor $F:\C^o \to \Cat$ to the category $\Cat$ of small
categories, let $\Tot(F)$ be the category of pairs $\langle c,s
\rangle$ of an object $c \in \C$ and an object $s \in F(C)$, with
morphisms from $\langle c,s \rangle$ to $\langle c',s' \rangle$
given by a morphism $f:c \to c'$ and a morphism $s \to
F(f)(s)$. Then the forgetful functor $\pi:\Tot(F) \to \C$ is a
fibration, with fibers $\pi_c \cong F(c)$, $c \in \C$. If $F$ is a
functor to sets, so that for any $c \in \C$, $F(c)$ is a discrete
category, then the fibration $\pi$ is discrete.

Conversely, for any fibration $\pi:\C' \to \C$ with of small
categories, and any object $c \in \C$, let $\G(\pi)(c)$ be the
category of cartesian functors $\C/c \to \C'$. Then $\G(\pi)(c)$ is
contravariantly functorial in $c$ and gives a functor $\G(\pi):\C^o
\to \Cat$. The two constructions are inverse, in the sense that we
have a natural cartesian equivalence $\Tot(\G(\pi)) \cong \C'$ for
any fibration $\pi':\C' \to \C$, and a natural pointwise equivalence
$F \to \G(\Tot(F))$ for any $F:\C^o \to \Cat$. In particular, we
have equivalences
$$
\pi_c \cong \G(\pi)(c), \qquad c \in \C.
$$
These equivalences of categories are not isomorphisms, so that the
fibers $\pi_c$ themselves do not form a functor from $\C^o$ to
$\Cat$ -- they only form a pseudofunctor in the sense of \cite{SGA}
(we do have a transition functor $f^*:\pi_{c'} \to \pi_c$ for any
morphism $f:c \to c'$ in $\C$, but this is compatible with
compositions only up to a canonical isomorphism). Nevertheless, for
all practical purposes, a fibered category over $\C$ is a convenient
axiomatization of the notion of a family of categories
contravariantly indexed by $\C$.

For any fibration $\pi:\C' \to \C$ of small categories, and any
functor $\gamma:\C_1 \to \C$ from a small category $\C_1$, we can
define a category $\gamma^*\C'$ and a functor $\pi_1:\gamma^*\C' \to
\C_1$ by taking the cartesian square
\begin{equation}\label{sq}
\begin{CD}
\gamma^*\C' @>{\gamma'}>> \C'\\
@V{\pi_1}VV @VV{\pi}V\\
\C_1 @>{\gamma}>> \C
\end{CD}
\end{equation}
in $\Cat$. Then $\pi_1$ is also a fibration, called the {\em induced
  fibration}. The corresponding pseudofunctor $\G(\pi_1):\C_1^o \to
\Cat$ is the composition of the functor $\gamma$ and $\G(\pi)$.

For covariantly indexed families, one uses the dual notion of a
cofibration: a morphism $f$ is cocartesian with respect to a functor
$\pi$ if it is cartesian with respect to $\pi^o$, a functor $\pi$ is
a cofibration if $\pi^o$ is a fibration, a functor $F:\C' \to \C''$
between two cofibrations is cocartesian if $F^o$ is cartesian, and a
subcategory $\C'_0 \subset \C'$ is a subcofibration if $(\C'_0)^o
\subset (\C')^o$ is a subfibration. The Grothendieck construction
associates cofibrations over $\C$ to functors from $\C$ to
$\Cat$. We have the same notion of an induced cofibration. Functors
to $\Sets \subset \Cat$ correspond to discrete cofibrations; the
simplest example of such is the projection
\begin{equation}\label{rho.c}
\rho_c:c\backslash\C \to \C
\end{equation}
for some object $c \in \C$, where $c\backslash\C = (\C^o/c)^o$ is
the category of objects $c' \in \C$ equipped with a map $c \to c'$.

\subsection{Base change.}\label{bc.subs}

Assume given a cofibration $\pi:\C' \to \C$ of small categories and
a functor $\gamma:\C_1 \to \C$, and consider the cartesian square
\eqref{sq}. Then the isomorphism $\gamma^{'*} \circ \pi^* \cong
\pi_1^* \circ \gamma^*$ induces by adjunction a base change map
$$
L^\hdot \pi_{1!} \circ \gamma^{'*} \to \gamma^* \circ L^\hdot \pi_!.
$$
This map is an isomorphism (for a proof see e.g. \cite{ka0}). In
particular, for any object $c \in \C$, any ring $R$, and any $E \in
\Fun(\C',R)$, we have a natural identification
\begin{equation}\label{pi.fib}
L^\hdot\pi_!E(c) \cong H_\idot(\pi_c,E|_c),
\end{equation}
where $E|_c \in \Fun(\pi_c,R)$ is the restriction to the fiber
$\pi_c \subset \C'$.  If the cofibration $\pi$ is discrete, then
this shows that $L^i\pi_!E = 0$ for $i \geq 1$, and
$$
\pi_!E(c) = \bigoplus_{c' \in \pi_c}E(c').
$$
For example, for the discrete cofibration $\rho_c$ of \eqref{rho.c}
and the constant functor $R \in \Fun(c\backslash\C,R)$, we obtain an
identification
\begin{equation}\label{r.rho.c}
R_c \cong \rho_{c!}R \cong L^\hdot\rho_{c!}R,
\end{equation}
where $R_c \in \Fun(\C,R)$ is the representable functor
\eqref{r.c}. For fibrations, we have exactly the same statements
with left Kan extensions replaced by right Kan extensions, and sums
replaced by products.

Moreover, assume that $R$ is commutative, and assume given an object
$T \in \Fun((\C')^o,R)$ that inverts all maps $f$ in $\C'$
cocartesian with respect to $\pi$ -- that is, $T(f)$ is invertible
for any such map. Then we can define the relative tensor product
functor $- \otimes_{\pi} T:\Fun(\C',R) \to \Fun(\C,R)$ by setting
$$
(E \otimes_{\pi} T)(c) = E|_c \otimes_{\pi_c} T|_c
$$
for any $E \in \Fun(\C',R)$. This has individual derived functors
$\Tor_\idot^\pi(-,T)$ and the total derived functor $ - \lotimes_{\pi}
T$. For any $c \in \C$, we have
\begin{equation}\label{ten.fib}
(E \lotimes_{\pi} T)(c) \cong E|_c \lotimes_{\pi_c} T|_c.
\end{equation}
If $T(c)$ is a free $R$-module for any $c \in \C'$, then we also
have the relative version
$$
\hhom_\pi(T,-):\Fun(\C,R) \to \Fun(\C',R)
$$
of the functor \eqref{hhom}; it is exact and right-adjoint to $-
\otimes_\pi T$, resp. $- \lotimes_\pi T$. In the case $T = \Z$, we
have $E \lotimes_\pi \Z \cong L^\hdot\pi_!E$, and the isomorphism
\eqref{ten.fib} is the isomorphism \eqref{pi.fib}.

\subsection{Simplicial objects.}\label{del.subs}

As usual, we denote by $\Delta$ the category of finite non-empty
totally ordered sets, a.k.a. finite non-empty ordinals, and somewhat
unusually, we denote by $[n] \in \Delta$ the set with $n$ elements,
$n \geq 1$. A simplicial object in a category $\C$ is a functor from
$\Delta^o$ to $\C$; these form a category denoted $\Delta^o\C$. For
any ring $R$ and $E \in \Fun(\Delta^o,R) = \Delta^oR\amod$, we
denote by $C_\idot(E)$ the normalized chain complex of the
simplicial $R$-module $E$. The homology of the complex $C_\idot(E)$
is canonically identified with the homology $H_\idot(\Delta^o,E)$ of
the category $\Delta^o$ with coefficients in $E$. Even stronger,
sending $E$ to $C_\idot(E)$ gives the {\em Dold-Kan equivalence}
$$
\Nn:\Fun(\Delta^o,R) \to C_{\geq 0}(R)
$$
between the category $\Fun(\Delta^o,R)$ and the category $C_{\geq
  0}(R)$ of complexes of $R$-modules concentrated in non-negative
homological degrees. The inverse equivalence is given by the
denormalization functor $\Dd:C_{\geq 0}(R) \to \Fun(\Delta^o,R)$
right-adjoint to $\Nn$.

For any simplicial set $X$, its homology $H_\idot(X,R)$ with
coefficients in a ring $R$ is the homology of the chain complex
$$
C_\idot(X,R) = C_\idot(R[X]),
$$
where $R[X] \in \Fun(\Delta^o,R)$ is given by $R[X]([n]) =
R[X([n])]$, $[n] \in \Delta$. By adjunction, for any simplicial set
$X$ and any complex $E_\idot \in C_{\geq 0}(R)$, a map $C_\idot(X,R)
\to E_\idot$ gives rise to a map of simplicial sets
\begin{equation}\label{taut}
\begin{CD}
X @>>> R[X] @>>> \Dd(E_\idot),
\end{CD}
\end{equation}
where we treat simplicial $R$-modules $R[X]$ and $\Dd(E_\idot)$ as
simplicial sets. Conversely, every map of simplicial sets $X \to
\Dd(E_\idot)$ gives rise to a map $C_\idot(X,R) \to E_\idot$. In
particular, if we take $X=\Dd(E_\idot)$, we obtain the {\em assembly
  map}
\begin{equation}\label{asse}
C_\idot(\Dd(E_\idot),R) \to E_\idot.
\end{equation}
The constructions are mutually inverse: every map of complexes of
$R$-modules $C_\idot(X,R) \to E_\idot$ decomposes as
\begin{equation}\label{ass.fa}
\begin{CD}
C_\idot(X,R) @>>> C_\idot(\Dd(E_\idot),R) @>>> E_\idot,
\end{CD}
\end{equation}
where the first map is induced by the tautological map \eqref{taut},
and the second map is the assembly map \eqref{asse}.

Applying the Grothendieck construction to a simplicial set $X$, we
obtain a category $\Tot(X)$ with a discrete fibration $\pi:\Tot(X)
\to \Delta$. We then have a canonical identification
\begin{equation}\label{tot.eq}
H_\idot(\Tot(X)^o,R) \cong H_\idot(\Delta^o,\pi_!R) \cong
H_\idot(\Delta^o,R[X]),
\end{equation}
so that $H_\idot(X,R)$ is naturally identified with the homology of
the small category $\Tot(X)^o$ with coefficients in the constant
functor $R$.

The {\em nerve} of a small category $\C$ is the simplicial set
$N(\C) \in \Delta^o\Sets$ such that for any $[n] \in \Delta$,
$N(\C)([n])$ is the set of functors from the ordinal $[n]$ to
$\C$. Explicitly, elements in $N(\C)([n])$ are diagrams
\begin{equation}\label{n.eq}
\begin{CD}
c_1 @>>> \dots @>>> c_n
\end{CD}
\end{equation}
in $\C$. We denote by $\N(\C) = \Tot(N(\C))$ the corresponding
fibered category over $\Delta$. Then by definition, objects of
$\N(\C)$ are diagrams \eqref{n.eq}, and sending such a diagram to
$c_n$ gives a functor
\begin{equation}\label{q.eq}
q:\N(\C) \to \C.
\end{equation}
Say that a map $f:[n] \to [m]$ in $\Delta$ is {\em special} if it
identifies $[n]$ with a terminal segment of the ordinal $[m]$. For
any fibration $\pi:\C' \to \Delta$, say that a map $f$ in $\C'$ is
{\em special} if it is cartesian with respect to $\pi$ and $\pi(f)$
is special in $\Delta$, and say that a functor $F:\C' \to \E$ to
some category $E$ is {\em special} if it $F(f)$ is invertible for
any special map $f$ in $\C'$. Then the functor $q$ of \eqref{q.eq}
is special, and any special functor factors uniquely through $q$. In
particular, $\Fun(\C,R)$ is naturally identified the full
subcategory in $\Fun(\N(\C),R)$ spanned by special
functors. Moreover, on the level of derived categories, say that $E
\in \D(\C',R)$ is {\em special} if so is $\D(E):\C' \to \D(R)$, and
denote by $\D_{sp}(\C',R) \subset \D(\C',R)$ the full subcategory
spanned by special objects. Then the pullback functor
\begin{equation}\label{q.st}
q^*:\D(\C,R) \to \D(\N(\C),R)
\end{equation}
induces an equivalence between $\D(\C,R)$ and
$\D_{sp}(\N(\C),R)$. In particular, we have a natural isomorphism
\begin{equation}\label{N.eq}
H_\idot(\C,R) \cong H_\idot(\N(\C),R),
\end{equation}
and by \eqref{tot.eq}, the right-hand side is also canonically
identified with the homology $H_\idot(N(\C),R)$ of the simplicial
set $N(\C)$.

The {\em geometric realization functor} $X \mapsto |X|$ is a functor
from $\Delta^o\Sets$ to the category $\Top$ of topological
spaces. For any simplicial set $X$ and any ring $R$, the homology
$H_\idot(X,R)$ is naturally identified with the homology
$H_\idot(|X|,R)$ of its realization, and the isomorphism
\eqref{N.eq} can also be deduced from the following geometric fact:
for any simplicial set $X$, we have a natural homotopy equivalence
$$
|N(\Tot(X))| \cong |X|.
$$
Even stronger, the geometric realization functor extends to a
functor from $\Delta^o\Top$ to $\Top$, and for any small category
$\C$ equipped with a fibration $\pi:\C \to \Delta$, we have a
natural homotopy equivalence
\begin{equation}\label{N.G}
|N(\C)| \cong ||N(\G(\pi))||,
\end{equation}
where $N(\G(\pi)):\Delta^o \to \Delta^o\Sets$ is the natural
bisimplicial set corresponding to $\pi$, and $||-||$ in the
right-hand side stands for the geometric realization of its
pointwise geometric realization.  Geometric realization commutes
with products by the well-known Milnor Theorem, so that in
particular, \eqref{N.G} implies that for any self-product $\C
\times_\Delta \dots \times_\Delta \C$, we have a natural homotopy
equivalence
\begin{equation}\label{self.eq}
|N(\C \times_\Delta \dots \times_\Delta \C)| \cong |N(\C)| \times
\dots \times |N(\C)|.
\end{equation}

\subsection{$2$-categories.}\label{2.subs}

We will also need to work with $2$-categories, and for this, the
language of nerves is very convenient, since the nerve of a
$2$-category can be converted into a $1$-category by the Grothendieck
construction.

Namely, recall that a $2$-category $\C$ is given by a class of
objects $c \in \C$, a collection of morphism categories $\C(c,c')$,
$c,c' \in \C$, a collection of identity objects $\id_c \in \C(c,c)$
for any $c \in \C$, and a collection of composition functors
\begin{equation}\label{m.cc}
m_{c,c',c''}:\C(c,c') \times \C(c',c'') \to \C(c,c''), \qquad
c,c',c'' \in \C
\end{equation}
equiped with associativity and unitality isomorphisms, subject to
standard higher contraints. A $1$-category is then a $2$-category
$\C$ with discrete $\C(c,c')$, $c,c' \in \C$. For any $2$-category
$\C$ and any $[n] \in \Delta$, one can consider the category
$$
N(\C)_n = \coprod_{c_1,\dots,c_n \in \C}\C(c_1,c_2) \times \dots
\times \dots \C(c_{n-1},c_n).
$$
If $\C$ is a small $1$-category, then $N(\C)_n = N(\C)([n])$ is the
value of the nerve $N(\C) \in \Delta^o\Sets$ at $[n] \in \Delta$,
and the structure maps of the functor $N(\C):\Delta^o \to \Sets$ are
induced by the composition and unity maps in $\C$. In the general
case, the composition and unity functors turn $N(\C)$ into a
pseudofunctor from $\Delta^o$ to $\Cat$. We let
$$
\N(\C) = \Tot(N(\C))
$$
be the corresponding fibered category over $\Delta$, and call it the
{\em nerve} of the $2$-category $\C$.

The associativity and unitality isomorphisms in $\C$ give rise to
the compatibility isomorphisms of the pseudofunctor $N(\C)$, so that
they are encoded by the fibration $\N(\C) \to \Delta$. One can in
fact use this to give an alternative definition of a $2$-category,
see e.g. \cite{trace}, but we will not need this. However, it is
useful to note what happens to functors. A {\em $2$-functor} $F:\C
\to \C'$ between $2$-categories $\C$, $\C'$ is given by a map $F$
between their classes of objects, a collection of functors
\begin{equation}\label{F.cc}
F(c,c'):\C(c,c') \to \C'(F(c),F(c')), \qquad c,c' \in \C,
\end{equation}
and a collection of isomorphisms $F(c,c)(\id_c) \cong
\id_{F(c)}$, $c \in \C$, and
$$
m_{F(c),F(c'')} \circ (F(c,c') \times F(c',c'')) \cong F(c,c'')
\circ m_{c,c',c''}, \qquad c,c',c'' \in \C,
$$
again subject to standard higher constraints. Such a $2$-functor
gives rise to a functor $\N(F):\N(\C) \to
\N(\C')$ cartesian over $\Delta$, and the correspondence between
$2$-functors and cartesian functors is one-to-one.

The category $\Cat$ is a $2$-category in a natural way, and the
Grothendieck construction generalizes directly to $2$-functors from
a $2$-category $\C$ to $\Cat$. Namely, say that a cofibration
$\pi:\C' \to \N(\C)$ is {\em special} if for any special morphism
$f:c \to c'$ in $\N(\C)$, the transition functor $f_1:\pi_c \to
\pi_{c'}$ is an equivalence. Then $2$-functors $F:\C \to \Cat$
correspond to special cofibrations $\Tot(F) \to \N(\C)$, and the
correspondence is again one-to-one. If $\C$ is actually a
$1$-category, then a $2$-functor $F:\C \to \Cat$ is exactly the same
thing as a pseudofunctor $\overline{F}:\C \to \Cat$ in the sense of
the usual Grothendieck construction, and we have $\Tot(F) \cong
q^*\Tot(\overline{F})$, where $q$ is the functor of \eqref{q.eq}
(one easily checks that every special cofibration over $\N(\C)$
arises in this way).

The simplest example of a $2$-functor from a $2$-category $\C$ to
$\Cat$ is the functor $\C(c,-)$ represented by an object $c \in
\C$. We denote the corresponding special cofibration by
\begin{equation}\label{rho.c.2}
\wt{\rho}_c:\N(c\backslash\C) \to \N(\C).
\end{equation}
If $\C$ is a $1$-category, then $\wt{\rho}_c = q^*\rho_c$, where
$\rho_c$ is the discrete cofibration \eqref{rho.c}

\subsection{Homology of $2$-categories.}\label{ho.2.subs}

To define the derived category of functors from a small $2$-category
$\C$ to complexes of modules over a ring $R$, we use its nerve
$\N(\C)$, with its fibration $\pi:\N(\C) \to \Delta$ and the
associated notion of a special map and a special object.

\begin{defn}
For any ring $R$ and small $2$-category $\C$, the {\em derived
  category of functors from $\C$ to $R\amod$} is given by
$$
\D(\C,R) = \D_{sp}(\N(\C),R).
$$
\end{defn}

Recall that if $\C$ is a $1$-category, then $\D_{sp}(\N(\C),R)$ is
identified with $\D(\C,R)$ by the functor $q^*$ of \eqref{q.st}, so
that the notation is consistent. Since the truncation functors with
respect to the standard $t$-structure on $\D(\N(\C),R)$ send special
objects to special objects, this standard $t$-structure induces a
$t$-structure on $\D(\C,R) \subset \D(\N(\C),R)$ that we also call
standard. We denote its heart by $\Fun(\C,R) \subset \D(\C,R)$; it
is equivalent to the category of special functors from $\N(\C)$ to
$R\amod$. If $\C$ is a $1$-category, every special functor factors
uniquely through $q$ of \eqref{q.eq}, so that the notation is still
consistent.

\begin{lemma}\label{sp.lemma}
For any $2$-category $\C$, the embedding $\D(\C,R) \subset
\D(\N(\C),R)$ admits a left and a right-adjoint functors
$L^{sp},R^{sp}:\D(\N(\C),R) \to \D(\C,R)$. For any object $c \in \C$
with the correspoding object $n(c) \in N(\C)_1 \subset \N(\C)$, we
have
$$
L^{sp}R_{n(c)} \cong L^\hdot\wt{\rho}_{c!}R,
$$
where $\wt{\rho}_c$ is the special cofibration \eqref{rho.c.2}, and
$R$ in the right-hand side is the constant functor.
\end{lemma}

\proof{} Say that a map $f$ in $\D(\N(\C))$ is {\em co-special} if
$\pi(f):[n] \to [n']$ sends the initial object of the ordinal $[n]$
to the initial object of the ordinal $[n']$. Then as in the proof of
\cite[Lemma 4.8]{mackey}, it is elementary to check that special and
co-special maps in $\N(C)$ form a complementary pair in the sense of
\cite[Definition 4.3]{mackey}, and then the adjoint functor $L^{sp}$
is provided by \cite[Lemma 4.6]{mackey}. Moreover, $L^{sp} \circ
L^{sp} \cong L^{sp}$, and $L^{sp}$ is an idempotent comonad on
$\D(\N(\C),R)$, with algebras over this comonad being exactly the
objects of $\D(\C,R)$. Moreover, by construction of \cite[Lemma
  4.6]{mackey}, $L^{sp}:\D(\N(C),R) \to \D(\N(\C),R)$ has a
right-adjoint functor $R^{sp}:\D(\N(\C),R) \to \D(\N(\C),R)$. By
adjunction, $R^{sp}$ is an idempotent monad, algebras over this
monad are objects in $\D(\C,R)$, and $R^{sp}$ factors through the
desired functor $\D(\N(\C,R)) \to \D(\C,R)$ right-adjoint to the
embedding $\D(\C,R) \subset \D(\N(\C),R)$. Finally, the last claim
immediately follows by the same argument as in the proof of
\cite[Theorem 4.2]{mackey}.
\endproof

For any $2$-functor $F:\C \to \C'$ between small $2$-categories, the
corresponding functor $\N(F)$ sends special maps to special maps, so
that we have a pullback functor
$$
F^* = \N(F)^*:\D(\C',R) \to \D(\C,R).
$$
By Lemma~\ref{sp.lemma}, $F^*$ has a left and a right-adjoint
functor $F_!$, $F_*$, given by
$$
F_! = L^{sp} \circ L^\hdot \N(F)_!, \qquad F_* = R^{sp} \circ
R^\hdot \N(F)_*.
$$
For any object $c \in \C$, we denote
\begin{equation}\label{r.c.2}
R_c = L^{sp}R_{n(c)} \cong L^\hdot\wt{\rho}_{c!}R \in \D(\C,R).
\end{equation}
If $\C$ is a $1$-category, then this is consistent with \eqref{r.c}
by \eqref{r.rho.c}. In the general case, by base change, we have a
natural identification
\begin{equation}\label{r.c.2.e}
R_c(c') \cong H_\idot(\C(c,c'),R)
\end{equation}
for any $c' \in \C$, an analog of \eqref{r.c}. Moreover, by
adjunction, we have a natural isomorphism
\begin{equation}\label{r.c.re.2}
E(c) \cong \Hom(R_c,E)
\end{equation}
for any $E \in \D(\C,R)$, a generalization of \eqref{r.c.re}.

\subsection{Finite sets.}\label{fin.subs}

The first example of a $2$-category that we will need is the
following. Denote by $\Gamma$ the category of finite sets. Then
objects of the $2$-category $\Q\Gamma$ are finite sets $S \in
\Gamma$, and for any two $S_1,S_2 \in \Gamma$, the category
$\Q\Gamma(S_1,S_2)$ is the groupoid of diagrams
\begin{equation}\label{domik}
\begin{CD}
S_1 @<{l}<< S @>{r}>> S_2
\end{CD}
\end{equation}
in $\Gamma$ and isomorphisms between them. The composition functors
\eqref{m.cc} are obtained by taking fibered products.

We can also define a smaller $2$-category $\Gamma_+ \subset
\Q\Gamma$ by keeping the same objects and requiring that
$\Gamma_+(S_1,S_2)$ consists of diagrams \eqref{domik} with
injective map $l$. Then since such diagrams have no non-trivial
automorphisms, $\Gamma_+$ is actually a $1$-category. It is
equivalent to the category of finite pointed sets. The equivalence
sends a set $S$ with a disntiguished element $o \in S$ to the
complement $\overline{S} = S \setminus \{o\}$, and a map $f:S \to
S'$ goes to the diagram
$$
\begin{CD}
\overline{S} @<{i}<< f^{-1}(\overline{S}') @>{f}>> \overline{S}',
\end{CD}
$$
where $i:f^{-1}(\overline{S}') \to \overline{S}$ is the natural
embedding. For any $n \geq 0$, we denote by $[n]_+ \in \Gamma_+$ the
set with $n$ non-distinguished elements (and one distinguished
element $o$). In particular, $[0]_+ = \{o\}$ is the set with the
single element $o$.

To construct $2$-functors from $\Q\Gamma$ to $\Cat$, recall that for
any category $\C$, the {\em wreath product} $\C\wr\Gamma$ is the
category of pairs $\langle S,\{c_s\}\rangle$ of a set $S \in \Gamma$
and a collection of objects $c_s \in \C$ indexed by elements $s \in
S$. Morphisms from $\langle S,\{c_s\} \rangle$ to $\langle
S',\{c'_s\}\rangle$ are given by a morphism $f:S \to S'$ and a
collection of morphisms $c_s \to c'_{f(s)}$, $s \in S$. Then the
forgetful functor $\rho:\C\wr\Gamma \to \Gamma$ is a fibration whose
fiber over $S \in \Gamma$ is the product $\C^S$ of copies of the
category $\C$ numbered by elements $s \in S$, and whose transition
functor $f^*:\C^{S_1} \to \C^{S_2}$ associated to a map $f:S_1 \to
S_2$ is the natural pullback functor.

Assume that the category $\C$ has finite coproducts (including the
coproduct of an empty collection of objects, namely, the initial
object $0 \in \C$). Then all the transition functors $f^*$ of the
fibration $\rho$ have left-adjoint functors $f_!$, so that $\rho$ is
also a cofibration. Moreover, for any diagram \eqref{domik} in
$\Gamma$, we have a natural functor
\begin{equation}\label{lr.eq}
r_! \circ l^*:\C^{S_1} \to \C^{S_2}.
\end{equation}
This defines a $2$-functor $\Vect(\C):\Q\Gamma \to \Cat$ -- for any
finite set $S \in \Gamma$, we let $\Vect(\C)(S) = \C^S$, and for any
$S_1,S_2 \in \Gamma$, the functor $\Vect(\C)(S_1,S_2)$ of
\eqref{F.cc} sends a diagram \eqref{domik} to the functor induced by
\eqref{lr.eq}. Moreover, for any subcategory $w(\C) \subset \C$ with
the same objects as $\C$ and containing all isomorphisms, the
collection of subcategories $\Vect(w(\C))(S) = w(\C)^S \subset
\C^S$ defines a subfunctor $\Vect(w(\C)) \subset \Vect(\C)$.

Restricting the $2$-functor $\Vect(\C)$ to the subcategory $\Gamma_+
\subset \Q\Gamma$ and applying the Grothendieck construction, we
obtain a cofibration over $\Gamma_+$ that we denote by
$\rho_+:(\C\wr\Gamma)_+ \to \Gamma_+$. For any subcategory $w(\C)$
with the same objects an containng all isomorphisms, we can do the
same procedure with the subfunctor $\Vect(w(\C)) \subset \Vect(\C)$;
this gives a subcofibration $(w(\C)\wr\Gamma)_+ \subset
(\C\wr\Gamma)_+$, and in particular, $\rho_+$ restricts to a
cofibration
\begin{equation}\label{wr}
\rho_+:(w(\C)\wr\Gamma)_+ \to \Gamma_+.
\end{equation}
Explicitly, the fiber of the cofibration $\rho_+$ over a pointed set
$S \in \Gamma_+$ is identified with $w(\C)^{\overline{S}}$, where
$\overline{S} \subset S$ is the complement to the distiguished
element. The transition functor corresponding to a pointed map $f:S
\to S'$ sends a collection $\{c_s\} \in
\overline{\C}^{\overline{S}}$, $s \in \overline{S}$ to the
collection $c'_{s'}$, $s' \in \overline{S}'$ given by
\begin{equation}\label{cofi}
c'_{s'} = \bigoplus_{s \in f^{-1}(s')}c_s,
\end{equation}
where $\oplus$ stands for the coproduct in the category $\C$.

\subsection{Matrices and vectors.}\label{mat.subs}

Now more generally, assume that we are given a small category $\C_0$
with finite coproducts and initial object, and moreover, $\C_0$ is a
unital monoidal category, with a unit object $1 \in \C_0$ and the
tensor product functor $- \otimes -$ that preserves finite
coproducts in each variable. Then we can define a $2$-category
$\Mat(\C_0)$ in the following way:
\begin{enumerate}
\item objects of $\Mat(\C_0)$ are finite sets $S \in \Gamma$,
\item for any $S_1,S_2 \in \Gamma$, $\Mat(\C_0)(S_1,S_2) \subset
  \C^{S_1 \times S_2}$ is the groupoid of isomorphisms of the
  category $\C^{S_1 \times S_2}$,
\item for any $S \in \Gamma$, $\id_S \in \Mat(\C_0)(S,S)$ is given by
  $\id_S = \delta_!(p^*(1))$, where $p:S \to \ppt$ is the projection
  to the point, and $\delta:S \to S \times S$ is the diagonal
  embedding, and
\item for any $S_1,S_2,S_2 \in \Gamma$, the composition functor
  $m_{S_1,S_2,S_3}$ of \eqref{m.cc} is given by
$$
m_{S_1,S_2,S_3} = p_{2!} \circ \delta_2^*,
$$
where $p_2:S_1 \times S_2 \times S_3 \to S_1 \times S_3$ is the
product $p_2 = \id \times p \times \id$, and analogously, $\delta_2
= \id \times \delta \times \id$.
\end{enumerate}
In other words, objects in $\Mat(\C_0)(S_1,S_2)$ are matrices of
objects in $\C$ indexed by $S_1 \times S_2$, and the identity object
and the composition functors are induces by those of $\C$ by the
usual matrix multiplication rules. The associativity and unitality
isomorphisms are also induced by those of $\C_0$. It is
straightforward to check that this indeed defines a $2$-category; to
simplify notation, we denote its nerve by
$$
\mat(\C_0) = \N(\Mat(\C_0)).
$$
Moreover, assume given another small category $\C$ with finite
coproducts, and assume that $\C$ is a unital right module category
over the unital monoidal category $\C_0$ -- that is, we have the
action functor
\begin{equation}\label{act}
- \otimes -:\C \otimes \C_0 \to \C,
\end{equation}
preserving finite coproducts in each variable and equipped with the
relevant unitality and asociativity isomorphism. Then we can define
a $2$-functor $\Vect(\C,\C_0)$ from $\Mat(\C_0)$ to $\Cat$ that
sends $S \in \Gamma$ to $\C^S$, and sends an object $M \in
\Mat(\C_0)(S_1,S_2)$ to the functor $\C^{S_1} \to \C^{S_2}$ induced
by \eqref{act} via the usual rule of matrix action on vectors. We
denote the corresponding special cofibration over $\mat(\C_0)$ by
$\vect(\C,\C_0)$.
Moreover, given a subcategory $w(\C) \subset \C$ with the same
objects and containing all the isomorphisms, we obtain a subfunctor
$\Vect(w(\C),\C_0) \subset \Vect(w(\C),\C_0)$ given by
$$
\Vect(w(\C),\C_0)(S) = w(\C)^S \subset \C^S = \Vect(w(\C),\C_0)(S).
$$
We denote the corresponding subcofibration by
$$
\vect(w(\C),\C_0) \subset \vect(\C,\C_0).
$$
If we take $\C_0=\Gamma$, and let $- \otimes -$ be the cartesian
product, then $\Mat(\C_0)$ is exactly the category $\Q\Gamma$ of
Subsection~\ref{fin.subs}. Moreover, any category $\C$ that has
finite coproducts is automatically a module category over $\Gamma$
with respect to the action functor
$$
c \otimes S = \bigoplus_{s \in S}c, \qquad c \in \C,S \in \Gamma,
$$
and we have $\Vect(\C,\Gamma) = \Vect(\C)$, $\Vect(w(\C),\Gamma) =
\Vect(w(\C))$. This example is universal in the following sense: for
any associative unital category $\C_0$ with finite coproducts, we
have a unique coproduct-preserving tensor functor $\Gamma \to \C_0$,
namely $S \mapsto 1 \otimes S$, so that we have a canonical
$2$-functor
\begin{equation}\label{ga.ma}
e:\Q\Gamma \to \Mat(\C_0).
\end{equation}
For any $\C_0$-module category $\C$ with finite coproducts, we have
a natural equivalence $e \circ \Vect(\C,\C_0) \cong \Vect(\C)$, and
similarly for $w(\C)$.

\subsection{The relative setting.}\label{rel.subs}

Finally, let us observe that the $2$-functors $\Vect(\C,\C_0)$,
$\Vect(w(\C),\C_0)$ can also be defined in the relative
situation. Namely, assume given a cofibration $\pi:\C \to \C'$ whose
fibers $\pi_c$, $c \in \C'$ have finite coproducts. Moreover, assume
that $\C$ is a module category over $\C_0$, and the action functor
\eqref{act} commutes with projections to $\C'$, thus induces
$\C_0$-module category structures on the fibers $\pi_c$ of the
cofibration $\pi$. Furthermore, assume that the induced action
functors on the fibers $\pi_c$ preserve finite coproducts in each
variable. Then we can define a natural $2$-functor
$\Vect(\C/\C',\C_0):\Mat(\C_0) \to \Cat$ by setting
\begin{equation}\label{vect.rel}
\Vect(\C/\C',\C_0)(S) = \C \times_{\C'} \dots \times_{\C'} \C
\end{equation}
where the terms in the product in the right-hand side are numbered
by elements of the finite set $S$. As in the absolute situation, the
categories $\Mat(\C_0)(S_1,S_2)$ act by the vector multiplication
rule. We denote by
$$
\vect(\C/\C',\C_0) \to \mat(\C_0)
$$
the special cofibration corresponding to the $2$-functor
$\Vect(\C/\C',\C_0)$, and we observe that the cofibration $\pi$
induces a natural cofibration
\begin{equation}\label{vect.fib}
\vect(\C/\C',\C_0) \to \C
\end{equation}
whose fiber over $c \in \C$ is naturally identified with
$\vect(\pi_c,\C_0)$. Moreover, if we have a subcategory $w(\C)
\subset \C$ with the same objects that contains all the
isomorphisms, and $w(\C) \subset \C$ is a subcofibration, then we
can let
$$
\Vect(w(\C)/\C',\C_0)(S) = w(\C) \times_{\C'} \dots \times_{\C'} w(\C)
\subset \Vect(\C/\C',\C_0)(S)
$$
for any finite set $S \in \Gamma$, and this gives a subfunctor
$\Vect(w(\C)/\C',\C_0) \subset \Vect(\C/\C',\C_0)$ and a subcofibration
$\vect(w(\C)/\C',\C_0) \subset \vect(\C/\C',\C_0)$. The cofibration
\eqref{vect.fib} then induces a cofibration
\begin{equation}\label{vect.fib.w}
\vect(w(\C)/\C',\C_0) \to \C
\end{equation}
whose fibers are identified with $\vect(w(\pi_c),\C_0)$, $c \in
\C$. As in the absolute case, in the case $\C_0=\Gamma$, we simplify
notation by setting $\vect(w(\C)/\C') = \vect(w(\C)/\C',\Gamma)$,
and we denote by
\begin{equation}\label{wr.rel}
((w(\C)/\C')\wr\Gamma)_+ \to \Gamma_+
\end{equation}
the induced cofibration over $\Gamma_+ \subset \Q\Gamma$.

Analogously, if $\pi:\C \to \C'$ is a fibration, then the same
constructions go through, except that $w(\C) \subset \C$ has to be a
subfibration, and the functors \eqref{vect.fib}, \eqref{vect.fib.w}
are also fibrations, with the same identification of the fibers.

\section{Statements.}\label{sta.sec}

\subsection{Generalities on $K$-theory.}

To fix notations and terminology, let us summarize very briefly the
definitions of algebraic $K$-theory groups.

First assume given a ring $k$, let $k\amod^{fp} \subset k\amod$ be
the category of finitely generated projective $k$-modules, and let
$BGL(k) \subset k\amod^{fp}$ be the groupoid of finitely generated
projective $k$-modules and their isomorphisms. Explicitly, we
have
$$
BGL(k) \cong \coprod_{P \in k\amod^{fp}}[\ppt/\Aut(P)],
$$
where the sum is over all isomorphism classes of finitely generated
projective $k$-modules, $\Aut(P)$ is the automorphism group of the
module $P$, and for any group $G$, $[\ppt/G]$ stands for the
groupoid with one object with automorphism group $G$. The category
$k\amod^{fp}$ is additive. In particular, it has finite
coproducts. Since $BGL(k) \subset k\amod^{fp}$ contains all objects
and all the isomorphisms, we have the cofibration
$$
\rho_+:(BGL(k)\wr\Gamma)_+ \to \Gamma_+
$$
of \eqref{wr}. Its fiber $(\rho_+)_{[1]_+}$ over the set $[1]_+
\in\Gamma_+$ is $BGL(k)$, and the fiber $(\rho_+)_S$ over a general
$S \in \Gamma_+$ is the product $BGL(k)^{\overline{S}}$. Applying
the Grothendieck construction and taking the geometric realization
of the nerve, we obtain a functor
$$
|N(\G(\rho_+))|:\Gamma_+ \to \Top
$$
from $\Gamma_+$ to the category $\Top$ of topological spaces, or in
other terminology, a $\Gamma$-space.  Then \eqref{cofi} immediately
shows that this $\Gamma$-space is special in the sense of the Segal
machine \cite{segal}, thus gives rise to a spectrum $\K(k)$. The
algebraic $K$-groups $K_\idot(k) = \pi_\idot\K(k)$ are by definition
the homotopy groups of this spectrum.

For a more general $K$-theory setup, assume given a small category
$\C$ with the subcategories $c(\C),w(\C) \subset \C$ of cofibrations
and weak equivalences, and assume that $\langle \C,c(\C),w(\C)
\rangle$ is a Waldhausen category. In particular, $\C$ has finite
coproducts and the initial object $0 \in \C$. Then one lets $E\C$ be
the category of pairs $\langle [n],\phi \rangle$ of an object $[n]
\in \Delta$ and a functor $\phi:[n] \to \C$, with morphisms from
$\langle [n],\phi \rangle$ to $\langle [n'],\phi' \rangle$ given by
a pair $\langle f,\alpha \rangle$ of a map $f:[n] \to [n']$ and a
map $\alpha:\phi' \circ f \to \phi$. Further, one lets $\wt{S\C}
\subset E\C$ be the full subcategory spanned by pairs $\langle [n],
\phi \rangle$ such that $\phi$ factors through $c(\C) \subset \C$
and sends the initial object $o \in [n]$ to $0 \in \C$. The
forgetful functor $s:\wt{S\C} \to \Delta$ sending $\langle [n],\phi
\rangle$ to $[n]$ is a fibration; explicitly, its fiber over $[n]
\in \Delta$ is the category of diagrams \eqref{n.eq} such that all
the maps are cofibrations, and $c_1=0$. Finally, one says that a map
$f$ in $\wt{S\C}$ is {\em admissible} if in its canonical
factorization $f = g \circ f'$ with $s(f) = s(f')$ and $f'$
cartesian with respect to $s$, the morphism $g$ pointwise lies in
$w(\C) \subset \C$. Then by definition, $S\C \subset \wt{S\C}$ is the
subcategory with the same objects and admissible maps between
them. This is again a fibered category over $\Delta$, with the
fibration $S\C \to \Delta$ induced by the forgetful functor $s$. The
$K$-groups $K_\idot(\C)$ are given by
$$
K_i(\C) = \pi_{i+1}(|N(S\C)|), \qquad i \geq 0.
$$
Moreover, since $\C$ has finite coproducts, the fibers of the
fibration $\wt{S\C} \to \Delta$ also have finite coproducts, and
since $S\C \subset \wt{S\C}$ contains all objects and all
isomorphisms, we can form the cofibration
\begin{equation}\label{rho.s}
\rho_+:((S\C/\Delta)\wr\Gamma)_+ \to \Gamma_+
\end{equation}
of \eqref{wr.rel}. Its fibers are the self-products $S\C
\times_\Delta \dots \times_\Delta S\C$. Then by \eqref{self.eq},
$$
|N(\G(\rho_+))|:\Gamma_+ \to \Top
$$
is a special $\Gamma$-space, so that $|N(S\C)|$ has a natural
infinite loop space structure and defines a connective spectrum. The
$K$-theory spectrum $\K(\C)$ is given by $\K(\C) = \Omega
|N(S\C)|$.

\begin{remark}
Our definition of the category $S\C$ differs from the usual one in
that the fibers of the fibration $s$ are opposite to what one gets
in the usual definition. This is harmless since passing to the
opposite category does not change the homotopy type of the nerve,
and this allows for a more succint definition.
\end{remark}

The main reason we have reproduced the $S$-construction instead of
using it as a black box is the following observation: the
construction works just as well in the relative setting. Namely, let
us say that a {\em family of Waldhausen categories} indexed by a
category $\C'$ is a category $\C$ equipped with a cofibration
$\pi:\C \to \C'$ with small fibers, and two subcofibrations
$s(\C),w(\C) \subset \C$ such that for any $c \in \C'$, the
subcategories
$$
c(\pi_c) = c(\C) \cap \pi_c \subset \pi_c, \qquad w(\pi_c) = w(\C)
\cap \pi_c \subset \pi_c
$$
in the fiber $\pi_c$ of the cofibration $\pi$ turn it into a
Waldhausen category. Then given such a family, one defines the
category $E\C$ exactly as in the absolute case, and one lets
$\wt{S(\C/\C')} \subset E\C$ be the full subcategory spanned by
$\wt{S\pi_c} \subset E\pi_c \subset E\C$, $c \in \C'$. Further, one
observes that the forgetful functor $s:\wt{S(\C/\C')} \to \Delta$ is
a fibration, and as in the absolute case, one let $S(\C/\C') \subset
\wt{S(\C/\C')}$ be the subcategory spanned by maps $f$ in whose
canonical factorization $f = g \circ f'$ with $s(f) = s(f')$ and
$f'$ cartesian with respect to $s$, the morphism $g$ pointwise lies
in $w(\C) \subset \C$. One then checks easily that the cofibration
$\pi$ induces a cofibration
$$
S(\C/\C') \to \C'
$$
whose fiber over $c \in \C'$ is naturally identified with
$S\pi_c$. This cofibration is obviously functorial in $\C'$: for any
functor $\gamma:\C'' \to \C'$ with the induced cofibration
$\gamma^*\C \to \C''$, we have $S(\gamma^*\C/\C'') \cong
\gamma^*S(\C/\C')$.

\subsection{The setup and the statement.}

Now assume given a commutative ring $k$, so that $k\amod^{fp}$ is a
monoidal category, and a Waldhausen category $\C$ that is additive
and $k$-linear, so that $\C$ is a module category over
$k\amod^{fp}$. Then all the fibers of the fibration $S\C \to \Delta$
are also module categories over $k\amod^{fp}$. To simplify notation,
denote
$$
\mat(k) = \mat(k\amod^{fp}), \quad \kk(\C,k) =
\vect(S\C/\Delta,k\amod^{fp}).
$$
More generally, assume given a family $\pi:\C \to \C'$ of Waldhausen
categories, and assume that all the fibers of the cofibration $\pi$
are additive and $k$-linear, and transition functors are additive
$k$-linear functors. Then $\C$ is a $k\amod^{fp}$-module category
over $\C$, and we can form the cofibration
$$
\kk(\C/\C',k) = \vect(S(\C/\C')/\Delta,k\amod^{fp}) \to \C' \times
\mat(k).
$$
Denote by
\begin{equation}\label{p.12}
\wt{\pi}:\kk(\C/\C',k) \to \C', \qquad \phi:\kk(\C/\C',k) \to \mat(k)
\end{equation}
its compositions with the projections to $\C'$ resp. $\mat(k)$. Then
the fiber of the cofibration $\wt{\pi}$ over $c \in \C'$ is naturally
idenitified with the category $\kk(\pi_c,k)$.

\begin{defn}\label{ada}
Let $R$ be the localization of $\Z$ in a set of primes.  A
commutative ring $k$ is {\em $R$-adapted} if $K_i(k) \otimes R = 0$
for $i \geq 1$, and $K_0(k) \otimes R \cong R$ as a ring.
\end{defn}

\begin{exa}
Let $k$ be a finite field of characteristic $\cchar(k)=p$, and let
$R = \Z_{(p)}$ be the localization of $\Z$ in the prime ideal $p\Z
\subset \Z$. Then $k$ is $R$-adapted by the famous theorem of
Quillen \cite{qui}.
\end{exa}

Assume given an $R$-adapted commutative ring $k$. Any additive map
$K_0(k) \to R$ induces a map of spectra
\begin{equation}\label{rk.sp}
\K(k) \to H(R),
\end{equation}
where $H(R)$ is the Eilenberg-Maclane spectrum corresponding to $R$,
so that fixing an isomorphism $K_0(k) \otimes R \cong R$ fixes a map
\eqref{rk.sp}. Do this, and for any $P \in k\amod^{fp}$, denote by
$\rk(P) \in R$ the image of its class $[P] \in K_0(k) \subset K_0(k)
\otimes R$ under the isomorphism we have fixed. Let $M(R)$ be the
category of free finitely generated $R$-modules, and let $T \in
\Fun(M(R)^o,R)$ be the functor sending a free $R$-module $M$ to $M^*
= \Hom_R(M,R)$. Equivalently, objects in $M(R)$ are finite sets $S$,
and morphisms from $S_1$ to $S_2$ are elements in the set $R[S_1
  \times S_2]$. In this description, sending $P \in k\amod^{fp}$ to
$\rk(P)$ defines a $2$-functor $\Rk:\Mat(k) \to M(R)$. By abuse of
notation, we denote
$$
\rk = q \circ \N(\Rk):\mat(k) \to \N(M(R)) \to M(R).
$$
Since the projection $\phi$ of \eqref{p.12} obviously inverts all
maps cocartesian with respect to the cofibration $\pi_1$, the
pullback $\phi^*\rk^*T \in \Fun(\kk(\C/\C',k),R)$ also inverts all
such maps. Therefore we are in the situation of
Subsection~\ref{bc.subs}, and we have a well-defined object
\begin{equation}\label{k.r}
K^R_\idot(\C/C',k) = \Z \lotimes_{\pi_1} \phi^*\rk^*T \in \D(\C',R),
\end{equation}
where $\Z$ on the left-hand side of the product is the constant
functor with value $\Z$. If $\C'=\ppt$ is the point category, we
simplify notation by letting $K^R_\idot(\C,k) =
K^R_\idot(\C/\ppt,k)$. The object $K^R_\idot(\C/\C',k)$ is clearly
functorial in $\C'$: for any functor $\gamma:\C'' \to \C'$, we have
a natural isomorphism
$$
\gamma^*K^R_\idot(\C/C',k) \cong K^R_\idot(\gamma^*\C/C'',k).
$$
In particular, the value of $K^R_\idot(\C/\C',k)$ at an object $c
\in \C'$ is naturally identified with $K^R_\idot(\pi_c,k)$.  Here
is, then, the main result of the paper.

\begin{theorem}\label{main}
Assume given a $k$-linear additive small Waldhausen category $\C$,
and a ring $R$ that is $k$-adapted in the sense of
Definition~\ref{ada}, and let $\K^R(\C,k)$ be the Eilenberg-Maclane
spectum associated to the complex $K^R_\idot(\C,k)$ of
\eqref{k.r}. Then there exists a natural map of spectra
$$
\K(\C) \to \K^R(\C,k)
$$
that induces an isomorphism of homology with coefficients in $R$.
\end{theorem}

Here as usual, we define ``homology with coefficients in $R$'' of a
spectrum $X$ by
$$
H_\idot(X,R) = \pi_\idot(X \wedge H(R)).
$$
If $R$ is the localization of $\Z$ in the set of primes $S$,
then by the standard spectral sequence argument, Theorem~\ref{main}
implies that $\nu$ becomes a homotopy equivalence after localizing
at the same set of primes $S$.

\section{Proofs.}\label{proof.sec}

\subsection{Additive functors.}

Before we prove Theorem~\ref{main}, we need a couple of technical
facts on the categories $\D(\Mat(k),R)$, $\D(M(R),R)$. Recall that
we have a natural $2$-functor $e:\Q\Gamma \to \Mat(k)$ of
\eqref{ga.ma}. Composing it with the natural embedding $\Gamma_+ \to
\Q\Gamma$, we obtain a $2$-functor
$$
i:\Gamma_+ \to \Mat(k).
$$
Composing further with the $2$-functor $\Rk:\Mat(k) \to M(R)$, we
obtain a functor
$$
\ii:\Gamma_+ \to M(R).
$$
Explicitly, $\ii$ sends a finite pointed set $S$ to its reduced
span
$$
\ii(S) = \overline{R[S]} = R[S]/R \cdot \{o\},
$$
where $o \in S$ is the distinguished element. The object $T \in
\Fun(M(R)^o,R)$ gives by pullback objects $\rk^{o*}T \in
\Fun(\mat(k)^o,R)$, $\ii^{o*}T \in \Fun(\Gamma_+^o,R)$. For any $E \in
\D(\Gamma_+,R)$, denote
\begin{equation}\label{h.ga}
H_\idot^\Gamma(E) = \Tor^{\Gamma_+}_\idot(E,\ii^*T).
\end{equation}
Say that an object $E \in \D(\Gamma_+,R)$ is {\em pointed} if
$E([0]_+)=0$, where $[0]_+ = \{o\} \in \Gamma_+$ is the pointed set
consisting of the distinguished element.

\begin{lemma}\label{loop}
\begin{enumerate}
\item For any two pointed objects $E_1,E_2 \in \D(\Gamma_+,R)$, we
  have $H_\idot^\Gamma(E_1 \lotimes E_2) = 0$.
\item Assume given a spectrum $X$ represented by a $\Gamma$-space
  $|X|:\Gamma_+ \to \Top$ special in the sense of Segal, and let
  $C_\idot(|X|,R) \in \D(\Gamma_+,R)$ be the object obtained by
  taking pointwise the singular chain homology complex
  $C_\idot(-,R)$. Then there exists a natural identification
$$
H_\idot^\Gamma(C_\idot(|X|,R)) \cong H_\idot(X,R).
$$
\end{enumerate}
\end{lemma}

\proof{} Although both claims are due to T. Pirashvili, in this form,
\thetag{i} is \cite[Lemma 2.3]{loop}, and its corollary \thetag{ii}
is \cite[Theorem 3.2]{loop}.
\endproof

The category $\Gamma_+$ has coproducts -- for any $S,S' \in
\Gamma_+$, their coproduct $S \vee S' \in \Gamma_+$ is the disjoint
union of $S$ and $S'$ with distinguished elements glued together. The
embedding $S \to S \vee S'$ admits a canonical retraction $p:S \vee
S' \to S$ identical on $S$ and sending the rest to the distiguished
element, and similarly, $S' \to S \vee S'$ has a canonical retraction
$p':S \vee S' \to S'$.

\begin{defn}
An object $E \in \D(\Gamma_+,R)$ is {\em additive} if for any $S,S'
\in \Gamma_+$, the natural map
\begin{equation}\label{add.eq}
E(S \vee S') \to E(S) \oplus E(S')
\end{equation}
induced by the retractions $p$, $p'$ is an isomorphism. An object
$E$ in the category $\D(\Mat(k),R)$ resp. $\D(M(R),R)$ is {\em
  additive} if so is $i^*E$ resp. $\ii^*E$.
\end{defn}

We denote by $\D_{add}(\Gamma_+,R)$, $\D_{add}(\Mat(k),R)$,
$\D_{add}(M(R),R)$ the full subcategories in $\D(\Gamma_+,R)$,
$\D(\Mat(k),R)$, $\D(M(R),R)$ spanned by additive objects. In fact,
$\D_{add}(\Gamma_+,R)$ is easily seen to be equivalent to
$\D(R)$. Indeed, $[0]_+ \in \Gamma_+$ is a retract of $[1]_+ \in
\Gamma_+$, so that we have a canonical direct sum decomposition
$$
R_1 \cong t \oplus R_0,
$$
where to simplify notation, we denote $R_n = R_{[n]_+} \in
\Fun(\Gamma_+,R)$, $n \geq 0$. Then for any pointed $E \in
\D(\Gamma_+,R)$, the adjunction map induces a map
\begin{equation}\label{t.M}
t \otimes M \to E,
\end{equation}
where $M = E([1])_+ \in \D(R)$. Any additive object is automatically
pointed, and the map \eqref{t.M} is an isomorphism if and only if
$E$ is additive. We actually have $t \otimes M \cong
\hhom(\ii^{o*}T,M) \cong \ii^*\hhom(T,M)$, so that the equivalence
$\D(R) \cong \D_{add}(\Gamma_+,R)$ is realized by the functor
$$
\ii^* \circ \hhom(T,-):\D(R) \overset{\sim}{\longrightarrow}
\D_{add}(\Gamma_+,R) \subset \D(\Gamma_+,R).
$$

\subsection{Adjunctions.}

By definition, $\ii^*$ sends $\D_{add}(M(R),R)$ and $i^*$ sends
$\D_{add}(\Mat(k),R)$ into $\D_{add}(\Gamma_+,R) \subset
\D(\Mat(k),R)$. It turns out that the same is true for their adjoint
functors $i_*$, $R^\hdot\ii_*$.

\begin{lemma}\label{i.ii.le}
\begin{enumerate}
\item For any additive $\overline{E} \in \D(\Gamma_+,R)$, the
  objects $R^\hdot\ii_*\overline{E} \in \D(M(R),R)$ and
  $i_*\overline{E} \in \D(\Mat(k),R)$ are additive.
\item For any additive $E \in \Fun(\Mat(k),R) \subset
  \D(\Mat(k),R)$, the adjunction map $E \to i_*i^*E$ is an
  isomorphism in homological degree $0$ with respect to the standard
  $t$-structure.
\end{enumerate}
\end{lemma}

\proof{} For the first claim, let $E = R^\hdot\ii_*\overline{E}$,
and note that we may assume that $\overline{E} = \ii^*\hhom(T,M)$
for some $M \in D(R)$. Therefore by adjunction, to check that
\eqref{add.eq} is an isomorphism, we need to check that the natural
map
$$
H_\idot^\Gamma(\ii^*R_{\ii(S)}) \oplus
H_\idot^\Gamma(\ii^*R_{\ii(S')}) \to H_\idot^\Gamma(R_{\ii(S \vee
  S')})
$$
is an isomorphism, where $H_\idot^\Gamma(-)$ is as in \eqref{h.ga},
$R_{\ii(S)}$, $R_{\ii(S')}$, $R_{\ii(S \vee S')}$ are the
representable functors \eqref{r.c}, and the map is induced by the
projections $p$, $p'$. For any $S,S_1 \in \Gamma_+$, we have
\begin{equation}\label{ii.r}
\ii^*R_{\ii(S)}(S_1) \cong R[\overline{S} \times \overline{S}_1].
\end{equation}
In particular, $\ii^*R_{\ii(S)}([0]_+) \cong R$ indepedently of $S$,
and the tautological projection $S \to [0]_+$ induces a functorial
map
$$
t:\ii^*R_{\ii([0]_+)} \to \ii^*R_{\ii(S)} \cong R
$$
in $\Fun(\Gamma_+,R)$ identical after evaluation at $[0]_+ \in
\Gamma_+$.  Moreover, we have
\begin{equation}\label{times}
\ii^*R_{\ii(S \vee S')} \cong \ii^*R_{\ii(S)} \otimes
\ii^*R_{\ii(S')},
\end{equation}
and under these identifications, the projections $p$, $p'$ induce
maps $\id \otimes t$ resp. $t \otimes \id$. Then to finish the
proof, in suffices to invoke Lemma~\ref{loop}~\thetag{i}.

For the object $i_*\overline{E}$, the argument is the same, but we
need to replace the representable functors $R_{\ii(S)}$,
$R_{\ii(S')}$, $R_{\ii(S \vee S')}$ by their $2$-category versions
of \eqref{r.c.2}, and \eqref{ii.r} becomes the isomorphism
$$
i^*R_{i(S)} \cong H_\idot(BGL(k)^{\overline{S} \times
  \overline{S_1}},R)
$$
provided by \eqref{r.c.2.e}. The corresponding version of
\eqref{times} then follows from the K\"unneth formula.

For the second claim, note that since we have already proved that
$i_*i^*\overline{E}$ is additive, it suffices to prove that the
natural map
$$
E([1]_+) \to i_*i^*E([1]_+)
$$
is an isomorphism in homological degree $0$. Again by
Lemma~\ref{loop}~\thetag{ii} and adjunction, this amount to checking
that the natural map
$$
H_0(\K(k),R) \to R
$$
induced by the rank map $\rk$ is an isomorphism. This follows from
Definition~\ref{ada} and Hurewicz Theorem.
\endproof

By definition, the functor $\rk^*$ also sends additive objects to
additive objects, but here the situation is even better.

\begin{lemma}\label{mat.m}
The functor $\rk_*:\D(\Mat(k),R) \to \D(M(R),R)$ sends additive
objects to additive objects, and $\rk^*$, $\rk_*$ induce mutually
inverse equivalences between $\D_{add}(\Mat(k),R)$ and
$\D_{add}(M(R),R)$.
\end{lemma}

\proof{} Assume for a moment that we know that for any additive $E
\in \D(\Mat(k),R)$, $\rk_*E$ is additive, and the adjunction map
$\rk^*\rk_* \to E$ is an isomorphism. Then for any additive $E \in
\D_{add}(M(R),R)$, the cone of the adjunction map $E \to
\rk_*\rk^*E$ is annihilated by $\rk^*$. Since the functor $\rk^*$ is
obviously conservative, $E \to \rk_*\rk^*E$ then must be an
isomorphism, and this would prove the claim.

It remains to prove that for any $E \in \D_{add}(\Mat(k),R)$,
$\rk_*E$ is additive, and the map $\rk^*\rk_*E \to E$ is an
isomorphism. Note that we have
$$
E \cong \dlim_{\overset{n}{\gets}}\tau_{\geq -n}E,
$$
where $\tau_{\geq -n}E$ is the truncation with respect to the
standard $t$-structure. If $E$ is additive, then all its truncations
are additive, and by adjunction, $\rk_*$ commutes with derived
inverse limits. Moreover, since derived inverse limit commutes with
finite sums, it preserves the additivity condition. Thus it suffices
to prove the statement under assumption that $E$ is bounded from
below with respect to the standard $t$-structure. Moreover, it
suffices to prove it separately in each homological degree $n$.

Since $\rk^*$ is obviously exact with respect to the standard
$t$-structure, $\rk_*$ is right-exact by adjunction, and the
statement is trivially true for $E \in \D^{\geq
  n+1}(\Mat(k),R)$. Therefore by induction, we may assume that the
statement is proved for $E \in \D^{\geq m+1}_{add}(\Mat(k),R)$ for
some $m$, and we need to prove it for $E \in \D^{\geq
  m}_{add}(\Mat(k),R)$. Let $\overline{E} = i^*E$. Since $E$ is
additive, $\overline{E}$ is also additive, so that $i_*\overline{E}$
is additive by Lemma~\ref{i.ii.le}~\thetag{i}.  The functor $i_*$ is
also right-exact with respect to the standard $t$-structures by
adjunction, and by Lemma~\ref{i.ii.le}~\thetag{ii}, the cone of the
adjunction map
$$
E \to i_*i^*E = i_*\overline{E}
$$
lies in $\D^{\geq m+1}_{add}(\Mat(k),R)$. Therefore it suffices to
prove the statement for $i_*\overline{E}$ instead of $E$. Since
$\rk_*i_*\overline{E} \cong R^\hdot\ii_*\overline{E}$ is additive by
Lemma~\ref{i.ii.le}~\thetag{i}, it suffices to prove that the
adjunction map
$$
\rk^*\ii_*\overline{E} \cong \rk^*\rk_*i_*\overline{E} \to i_*\overline{E}
$$
is an isomorphism. Moreover, since both sides are additive, it
suffices to prove it after evaluating at $i([1]_+)$. We may assume
that $\overline{E} = \hhom(\ii^*T,M)$ for some $M \in \D(R)$, so
that by adjunction, this is equivalent to proving that the natural
map
$$
H_\idot^\Gamma(i^*R_{i([1]_+)}) \to H_\idot^\Gamma(\ii^*R_{\ii([1]_+)})
$$
is an isomorphism. But as in the proof of Lemma~\ref{i.ii.le}, this
map is the map
$$
H_\idot^\Gamma(C_\idot(BGL^{\overline{S}},R))
\to H_\idot^\Gamma(R[\overline{S}])
$$
induced by the functor $\rk$, and by Lemma~\ref{loop}~\thetag{ii},
it is identified with the map of homology
$$
H_\idot(\K(k),R) \to H_\idot(H(R),R)
$$
induced by the map of spectra \eqref{rk.sp}. This map is an
isomorphism by Definition~\ref{ada}.
\endproof

\subsection{Proof of the theorem.}

We can now prove Theorem~\ref{main}. We begin by constructing the
map. To simplify notation, let $K = K^R_\idot(\C,k) \in \D(R)$, and
let
$$
E = L^\hdot\pi_{2!}R \in \D(\Mat(k),\Z) \subset \D(\mat(k),\Z).
$$
Then by the projection formula, we have a natural quasiisomorphism
$$
K \cong E \lotimes_{\mat(k)} \rk^{o*}T,
$$
so that by adjunction, we obtain a natural map
\begin{equation}\label{v.eq}
v:E \to \hhom(\rk^{o*}T,K).
\end{equation}
Restricting with respect to the $2$-functor $i:\Gamma_+ \to
\Mat(k)$, we obtain a map
\begin{equation}\label{bar.v}
\overline{v}:\overline{E} \to i^*\hhom(\rk^{o*}T,K) \cong
\hhom(\ii^{o*}T,K),
\end{equation}
where we denote $\overline{E} = i^*E$. Now note that over
$i(\N(\Gamma_+)) \subset \mat(k)$, the cofibration $\phi:\kk(\C,k)
\to \mat(k)$ restricts to the special cofibration corresponding to
the cofibration $\rho_+$ of \eqref{rho.s}. Therefore by base change,
we have $\overline{E} \cong L^\hdot\rho_{+!}R$. Then to compute
$\overline{E}$, we can apply the Grothendieck construction to the
cofibration $\rho_+$ and use base change; this shows that
$\overline{E} \in \D(\Gamma_+,R)$ can be represented by the homology
complex
$$
E_\idot = C_\idot(N(\G(\rho_+)),R).
$$
Choose a complex $\overline{K}_\idot$ representing $\hhom(\ii^*T,K)
\in \D(\Gamma_+,R)$ in such a way that the map $\overline{v}$ of
\eqref{bar.v} is represented by a map of complexes
$$
\overline{v}_\idot:E_\idot \to \overline{K}_\idot.
$$
Replacing $\overline{K}_\idot$ with its truncation if necessary, we
may assume that it is concentrated in non-negative homological
degrees. Applying the Dold-Kan equivalence pointwise, we obtain a
functor $\Dd(\overline{K}_\idot)$ from $\Gamma_+$ to simplicial
abelian groups. We can treat it as a functor to simplicial sets, and
take pointwise the tautological map \eqref{taut}; this results in a
map
\begin{equation}\label{nu.bar}
\overline{\nu}:N(\G(\rho_+)) \to \Dd(\overline{K}_\idot)
\end{equation}
of functors from $\Gamma_+$ to simplicial sets. Taking pointwise
geometric realization, we obtain a map of $\Gamma$-spaces, hence of
spectra. By definition, the $\Gamma$-space $|N(\G(\rho_+))|$
corresponds to the spectrum $\K(\C)$. Since $\overline{K}_\idot$
represents the additive object $\ii^*\hhom(T,K) \in \D(\Gamma_+,R)$,
the isomorphisms \eqref{add.eq} induce weak equivalences of
simplicial sets
$$
\Dd(\overline{K}_\idot)(S \vee S') \cong \Dd(\overline{K}_\idot)(S)
\times \Dd(\overline{K}_\idot)(S'),
$$
so that the $\Gamma$-space $|\Dd(\overline{K}_\idot)|$ is
special. It gives the Eilenberg-Maclane spectrum $\K$ corresponding
to $K \cong \overline{K}_\idot([1]_+) \in \D(R)$. Thus the map of
spectra induced by $\overline{\nu}$ of \eqref{nu.bar} reads as
\begin{equation}\label{thm.map}
\K(\C) \to \K.
\end{equation}
This is our map.

\medskip

To prove the theorem, we need to show that the map $\overline{\nu}$
induces an isomorphism on homology with coefficients in $R$. Let
$\overline{S} \in \D(\Gamma_+,R)$ be the object represented by the
chain complex $C_\idot(\Dd(\overline{K}_\idot),R)$. Then by
Lemma~\ref{loop}~\thetag{ii}, it suffices to prove that the map
\begin{equation}\label{e.s}
H_\idot^\Gamma(\overline{E}) \to H_\idot^\Gamma(\overline{S})
\end{equation}
induced by \eqref{nu.bar} is an isomorphism. Moreover, note that we
can apply the procedure above to the map $v$ of \eqref{v.eq} instead
of its restriction $\overline{v}$ of \eqref{bar.v}. This results
in a map of functors
$$
N(\G(\phi)) \to \Dd(K_\idot),
$$
where $K_\idot$ is a certain complex representing $\hhom(\rk^*T,K)
\in \D(\Mat(k),R)$. If we denote by $S \in \D(\Mat(k),R)$ the object
represented by $C_\idot(\Dd(K_\idot),R)$ and let
\begin{equation}\label{nu.eq}
\nu:E \to S
\end{equation}
be the map induced by the map $v$, then we have $S_0 \cong i^*S$,
$i^*\nu$ is the map induced by $\overline{\nu}$ of \eqref{nu.bar},
and \eqref{e.s} becomes the map
$$
H_\idot^\Gamma(i^*\nu):H_\idot^\Gamma(i^*E) \to H_\idot^\Gamma(i^*S).
$$
By adjunction and Lemma~\ref{i.ii.le}~\thetag{i}, it then suffices
to prove that for any additive $N \in \D(\Mat(k),R)$, the map
$$
\Hom(S,N) \to \Hom(E,N)
$$
induced by the map $\nu:E \to S$ is an isomorphism. By
Lemma~\ref{mat.m}, we may assume that $N \cong \rk^*\wt{N}$ for some
additive $\wt{N} \in \D(M(R),R)$, and by induction on degree, we may
further assume that $\wt{N}$ lies in a single homological
degree. But since $R$ is a localization of $\Z$, any additive
functor from $M(R)$ to $R$-modules is $R$-linear, thus of the form
$\hhom(T,M)$ for some $R$-module $M$. Thus we may assume $\wt{N} =
\hhom(T,M)$ for some $M \in \D(R)$. Again by adjunction, it then
suffices to prove that the map
$$
E \lotimes_{\mat(k)} \rk^{o*}T \to S \lotimes_{\mat(k)} \rk^{o*}T
$$
induced by the map $\nu$ of \eqref{nu.eq} is an isomorphism. But the
adjunction map $v$ of \eqref{v.eq} has the decomposition
\eqref{ass.fa} that reads as
$$
\begin{CD}
E @>{\nu}>> S @>{\kappa}>> \hhom(\rk^{o*}T,K),
\end{CD}
$$
where $\kappa$ is the assembly map \eqref{asse} for the complex
$K_\idot$. Thus to finish the proof, it suffices to check the
following.

\begin{lemma}
For any object $K \in \D(R)$ represented by a complex $K_\idot$ of
flat $R$-modules concentrated in non-negative homological degrees,
denote by $\wt{S} \in \D(M(R),R)$ the object represented by the
complex $C_\idot(\Dd(\hhom(T,K_\idot)),R)$, let $S = \rk^*\wt{S}$,
and let
$$
\rk^*\kappa:S \to \rk^*\hhom(T,K) \cong \hhom(\rk^{o*}T,K)
$$
be the pullback of the assembly map $\kappa:\wt{S} \to
\hhom(T,K)$. Then the map
$$
S \lotimes_{\mat(k)} \rk^{o*}T \to K
$$
adjoint to $\rk^*\kappa$ is an isomorphism in $\D(R)$.
\end{lemma}

\proof{} For any $M \in R\amod$, we can consider the functor
$\hhom(T,M)$ as a functor from $M(R)$ to sets, and we have the
assembly map
\begin{equation}\label{asso}
R[\hhom(T,M)] \to \hhom(T,M).
\end{equation}
If $M$ is finitely generated and free, then by definition, we have
$$
\begin{aligned}
R[\hhom(T,M)](M_1) &= R[\hhom(T,M)(M_1)] = R[\Hom(M_1^*,M)]\\
&\cong R[\Hom(M^*,M_1)]
\end{aligned}
$$
for any $M_1 \in M(R)$, so that $R[\hhom(T,M)] \cong R_{M^*}$ is a
representable functor. Therefore $\Tor_i^{M(R)}(R[\hhom(T,M)],T)$
vanishes for $i \geq 1$, and the map
$$
R[\hhom(T,M)] \lotimes_{M(R)} T \cong R[\hhom(T,M)] \otimes_{M(R)} T
\to M
$$
adjoint to the assembly map \eqref{asso} is an isomorphism. Since
$-\lotimes-$ commutes with filtered direct limits, the same is true
for an $R$-module $M$ that is only flat, not necessarily finitely
generated or free.

Moreover, consider the product $\Delta^o \times M(R) \to M(R)$, with
the projections $\tau:\Delta^o \times M(R) \to M(R)$,
$\tau':\Delta^o \times M(R) \to \Delta^o$. Then for any simplicial
pointwise flat $R$-module $M \in \Fun(\Delta^o,R)$, the map
\begin{equation}\label{a}
a:R[\hhom(\tau^*T,M)] \lotimes_{\tau'} \tau^*T \to M
\end{equation}
adjoint to the assembly map $R[\hhom(\tau^*T,M)] \to
\hhom(\tau^*T,M)$ is also an isomorphism. Apply this to $M =
\Dd(K_\idot)$, and note that we have
$$
K \cong L^\hdot\tau_!M, \qquad \wt{S} \cong
L^\hdot\tau_!R[\hhom(\tau^*T,M)],
$$
and the map $\wt{S} \lotimes_{M(R)} T \to K$ adjoint to the assembly
map $\kappa$ is exactly $L^\hdot\tau_!(a)$, where $a$ is the map
\eqref{a}. Therefore it is also an isomorphism.

To finish the proof, it remains to show that the natural map
$$
\wt{S} \lotimes_{M(R)} T \to \rk^*\wt{S} \lotimes_{\Mat(k)}
\rk^{o*}T = S \lotimes_{\Mat(k)} \rk^{o*}T
$$
is an isomorphism. By adjunction, it suffices to show that the
natural map
$$
\Hom(\wt{S},E) \to \Hom(\wt{S},\rk_*\rk^*E) \cong \Hom(S,\rk^*E)
$$
is an isomorphism for any additive $E \in \D(M(R),R)$, and this
immediately follows from Lemma~\ref{mat.m}.
\endproof

\bigskip

\noindent
{\sc
Steklov Math Institute, Algebraic Geometry section\\
Laboratory of Algebraic Geometry, NRU HSE\\
\mbox{}\hspace{30mm}and\\
IBS Center for Geometry and Physics, Pohang, Rep. of Korea}

\bigskip

\noindent
{\em E-mail address\/}: {\tt kaledin@mi.ras.ru}

\end{document}